\documentstyle[12pt,amsmath,amsfonts,amssymb]{article}
\newtheorem{theorem}{Theorem}[section]
\newtheorem{definition}[theorem]{Definition}
\newtheorem{corollary}[theorem]{Corollary}
\newtheorem{proposition}[theorem]{Proposition}
\newtheorem{remark}[theorem]{Remark}
\newtheorem{lemma}[theorem]{Lemma}

\newcommand {\Mc}      {{\mathcal M}}

\newcommand {\Ac}      {{\mathcal A}}

\newcommand {\Bc}      {{\mathcal B}}

\newcommand {\Pc}      {{\mathcal P}}
\newcommand {\Ec}      {{\mathcal E}}
\newcommand {\Qc}      {{\mathcal Q}}

\newcommand {\Ic}      {{\mathcal I}}
\newcommand {\WP}      {W^1_p(\RN)}
\newcommand {\WPO}     {W^1_p(\Omega)}

\newcommand {\WKP}     {W^k_p(\RN)}
\newcommand {\WKPO}    {W^k_p(\Omega)}
\newcommand {\WKQ}     {W^k_q(\RN)}
\newcommand {\WKQO}    {W^k_q(\Omega)}

\newcommand {\PK}      {\Pc_k}

\newcommand {\tC}      {\widetilde{C}}

\newcommand {\tQ}      {\widetilde{Q}}
\newcommand {\tI}      {\widetilde{I}}

\newcommand {\tq}      {\tilde{q}}
\newcommand {\tr}      {\tilde{r}}

\newcommand {\R}       {{\bf R}}
\newcommand {\N}       {{\bf N}}
\newcommand {\RN}      {\R^n}

\newcommand {\ve}      {\varepsilon}

\newcommand {\emp}     {\emptyset}
\newcommand {\cw}      {\curlywedge}

\newcommand {\GW}      {\widetilde{\Gamma}}
\newcommand {\intl}    {\int\limits}

\newcommand {\IM}      {\Ic_m}
\newcommand {\as}      {\alpha^*}
\newcommand {\ps}      {{p^*}}
\newcommand {\bz}      {\bar{z}}

\newcommand {\diam}    {\operatorname{diam}}
\newcommand {\dist}    {\operatorname{dist}}

\newcommand {\len}     {\operatorname{len}}
\newcommand {\lng}     {\operatorname{length}}

\newcommand {\bx}      {\hspace{10mm}$\Box$}

\newcommand {\nn}      {\nonumber}
\newcommand {\rf}[1]   {(\ref{#1})}   
\newcommand {\reff}[1] {\ref{#1}}     
\newcommand {\SECT}[2] {\section*{\centerline{\normalsize
{\bf #1}}} \setcounter{section}{#2}
\setcounter{theorem}{0}\setcounter{equation}{0}}

\newcommand{\lbl}[1]        {\label{#1}}            
\newcommand{\be}            {\begin{eqnarray}}
\newcommand{\bel}[1]        {\begin{eqnarray} \label{#1}}
\newcommand{\ee}            {\end{eqnarray}}
\begin{document}
\medskip
\centerline{\large{\bf On Sobolev extension domains in
$\RN$}}
\vspace*{10mm} \centerline{By  {\it Pavel Shvartsman}}
\vspace*{12 mm}
\renewcommand{\thefootnote}{ }
\footnotetext[1]{{\it\hspace{-6mm}Math Subject
Classification} 46E35}
\begin{abstract} We describe a class of Sobolev
$W^k_p$-extension domains $\Omega\subset \RN$ determined by
a certain inner subhyperbolic metric in $\Omega$. This
enables us to characterize finitely connected
Sobolev $W^1_p$-extension domains in ${\R}^{2}$ for each
$p>2$ .
\end{abstract}
\renewcommand{\thefootnote}{\arabic{footnote}}
\setcounter{footnote}{0}
\SECT{1. Introduction.}{1}
\indent
\par Let $\Omega$ be a domain in $\RN$. This paper is
devoted to the problem of extendability of functions from
the Sobolev space $\WKPO$ to functions from $\WKP$. We
recall that, given $k\in\N$ and $p\in[1,\infty]$, the
Sobolev space $\WKPO$, see e.g. Maz'ja \cite{M}, consists
of all functions $f\in L_{1,\,loc}(\Omega)$ whose
distributional partial derivatives on $\Omega$ of all
orders up to $k$ belong to $L_p(\Omega)$. $\WKPO$ is normed
by
$$
\|f\|_{\WKPO}:=
\sum\{\|D^\alpha f\|_{L_p(\Omega)}:|\alpha|\le k\}.
 $$
\par A domain $\Omega $ in ${\bf R}^{n}$ is said to be a
{\it Sobolev} $W_{p}^{k}$-{\it extension domain} if the
following isomorphism
$$ \WKPO=\WKP|_{\Omega } $$
holds. In other words, $\Omega $ is a Sobolev extension
domain (for the space $\WKP$) if every Sobolev function
$f\in \WKPO$ can be extended to a Sobolev
$W_{p}^{k}$-function $F$ defined on all of $\RN$. For
instance, Lipschitz domains (Calder\'{o}n \cite{C2},
$1<p<\infty$, Stein \cite{St}, $p=1,\infty$) in $\RN$ are
$W_{p}^{k}$-extension domains for every $p\in[1,\infty]$
and every $k\in\N$. Jones \cite{Jn} introduced a wider
class of $(\varepsilon ,\delta )$-domains and proved that
every $(\varepsilon ,\delta )$-domain is a Sobolev
$W_{p}^{k}$-extension domain in ${\bf R}^{n}$ for every
$k\ge 1$ and every $p\ge 1$. Burago and Maz'ya \cite{BM},
\cite{M}, Ch. 6, described extension domains for the space
$BV(\RN)$ of functions whose distributional derivatives of
the first order are finite Radon measures.
\par Our main result is the following
\begin{theorem} \lbl{MAIN-EXT} Let $n<p<\infty$ and let
$\Omega $ be a domain in ${\bf R}^{n}$. Suppose that there
exist constants  $C,\theta>0$ such that the following
condition is satisfied: for every $x,y\in\Omega$ such that
$\|x-y\|\le\theta$, there exists a rectifiable curve
$\gamma\subset\Omega$ joining $x$ to $y$ such that
\bel{F-IN} \int_{\gamma } \dist(z,\partial\Omega)
^{\frac{1-n}{p-1}}\,ds(z) \le\,
C\|x-y\|^{\frac{p-n}{p-1}}.\ee
Here $\partial\Omega$ denotes the boundary of $\Omega$ and
$ds$ denotes arc length measure.
\par Then $\Omega$ is a Sobolev $W_{q}^{k}$-extension
domain for every $k\ge 1$ and every $q>p^*$ where
$p^*\in(n,p)$ is a constant depending only on $n,p$ and
$C$.
\end{theorem}
\par For $k=1$ and $q>p$ this result has been proved by
Koskela \cite{K}.
\par Observe that this theorem is also known for the case
$p=\infty$ (with $\ps=q=\infty$). In that case every domain
$\Omega$ satisfying inequality \rf{F-IN} is
quasi-Euclidean, i.e., its inner metric is (locally)
equivalent to the Euclidean distance. This case was studied
by Whitney \cite{W3} who proved that every quasi-Euclidean
domain is a $W^k_\infty$-extension domain for every $k\ge
1$.
\par Our next result, Theorem \reff{RTWO-EXT}, relates to
description of Sobolev extension domains in $\R^2$. The
first result in this direction was obtained by
Gol'dstein, Latfullin and Vodop'janov \cite{GLV,GV1,GV2}
who proved that a finitely connected bounded planar domain
$\Omega $ is a Sobolev $W_{2}^{1}$-extension domain if and
only if $\Omega $ is an $(\ve,\delta)-$ domain in $\R^2$
for some $\ve,\delta>0$. Maz'ja \cite{M,MP} gave an example
of a simply connected domain $\Omega\subset\R^2$ such that
$\Omega$ is a $W^1_p$-extension domain for every $p\in[1,2)$,
while $\R^2\setminus \bar{\Omega}$ is a $W^1_p$-extension
domain for all $p>2$. However $\Omega$ is not an
$(\ve,\delta)-$domain for any $\ve$ and $\delta$.
\par Buckley and Koskela \cite{BKos} showed that if a
finitely connected bounded domain $\Omega\subset\R^2$ is a
Sobolev $W_{p}^{1}$-extension domain for some $p>2$, then
there exists a constant $C>0$ such that {\it for every}
$x,y\in\Omega$ there exists a rectifiable curve
$\gamma\subset\Omega$ satisfying inequality \rf{F-IN} (with
$n=2$). Combining this result with Theorem \reff{MAIN-EXT},
we obtain the following
\begin{theorem}\lbl{RTWO-EXT} Let $2<p<\infty$ and let
$\Omega$ be a finitely connected bounded planar domain.
Then $\Omega$ is a Sobolev $W_{p}^{1}$-extension domain if
and only if for some $C>0$ the following condition is
satisfied: for every $x,y\in\Omega$ there exists a
rectifiable curve $\gamma\subset\Omega$ joining $x$ to $y$
such that
$$
\int_{\gamma } \dist(z,\partial\Omega)
^{\frac{1}{1-p}}\,ds(z) \le\,
C\|x-y\|^{\frac{p-2}{p-1}}.
$$
\end{theorem}
\par We note that this result is also true for the case
$p=\infty$ and then the space $W^1_\infty$ can even be
replaced by $W^k_\infty$ for {\it arbitrary} $k\ge 1$. This
follows from the aforementioned theorem of Whitney
\cite{W3} combined with a result of Zobin \cite{Zob2} which
states that every finitely connected bounded planar
$W^k_\infty$-extension domain is quasi-Euclidean. Zobin
\cite{Zob1} also showed that for every $k\ge 2$ there
exists a bounded planar $W^k_\infty$-extension domain which
is not quasi-Euclidean.
\par Let us briefly indicate the main ideas of our approach
for the case $k=1$, i.e., for the Sobolev space $\WP$.
Recall that, when $p>n$, it follows from the Sobolev
embedding theorem that every function $f\in \WPO,p>n,$ can
be redefined, if necessary, on a subset of $\Omega $ of
Lebesgue measure zero so that it satisfies a {\it local}
H\"{o}lder condition of order $\alpha :=1-\frac{n}{p} $ on
$\Omega $: i.e., for every ball $B\subset \Omega $
\bel{S-Lip-OM}
|f(x)-f(y)|\le C(n,p)\|f\| _{\WPO}
\| x-y\| ^{1-%
\frac{n}{p}},~~~~x,y\in B.
\ee
\par We will identify each element of $\WPO$ with its
unique continuous representative. Thus we will be able to
restrict our attention to the case of continuous Sobolev
functions.
\par Following Buckley and Stanoyevitch \cite{BSt3}, given
$\alpha \in [0,1]$ and a rectifiable curve $\gamma \subset
\Omega $, we define the {\it subhyperbolic length} of
$\gamma $ by
$$
\len_{\alpha ,\Omega }(\gamma ) :=\int_{\gamma }
\dist(z,\partial\Omega)^{\alpha -1}\,ds(z).
$$
Then we let $d_{\alpha ,\Omega }$ denote the corresponding
{\it subhyperbolic metric} on $\Omega$ given, for each
$x,y\in \Omega $, by
\bel{DEF-D}
d_{\alpha ,\Omega }(x,y):=\inf_{\gamma }
\len_{\alpha ,\Omega }(\gamma )
\ee
where the infimum is taken over all rectifiable curves
$\gamma \subset \Omega $ joining $x$ to $y$.
\par The metric $d_{\alpha ,\Omega }$ was introduced and
studied by Gehring and Martio in \cite{GM}. See also
\cite{AHHL,L,BKos} for various further results using this
metric.  Note also that $\operatorname{len}_{0,\Omega }$
and $d_{0,\Omega }$ are the well-known {\it quasihyperbolic
length} and {\it quasihyperbolic distance}, and
$d_{1,\Omega }$ is the {\it inner (or geodesic) metric} on
$\Omega $.
\par The subhyperbolic metric $d_{\alpha ,\Omega }$ with
$\alpha =(p-n)/(p-1)$ arises naturally in the study of
Sobolev $W_{p}^{1}(\Omega )$-functions for $p>n$. In
particular, Buckley and Stanoyevitch \cite{BSt2} proved
that the local H\"{o}lder condition \rf{S-Lip-OM} is
equivalent to the following H\"{o}lder-type condition
\bel{H-LOC} |f(x)-f(y)|\le C(n,p)\| f\| _{\WPO}\{d_{\alpha
,\Omega }(x,y)^{1-\frac{1}{p}}+\| x-y\|
^{1-\frac{n}{p}}\},~~~~x,y\in \Omega , \ee
with $\alpha =(p-n)/(p-1)$.
\par In turn, since any extension $F\in \WP$ of $f$
satisfies the global H\"{o}lder condition
$$
|F(x)-F(y)|\le
C(n,p)\|F\|_{\WP}\|x-y\|^{1-\frac{n}{p}},~~~~~x,y\in \RN,
$$
we have
\bel{HTOM} |f(x)-f(y)|\le
C(n,p)\|F\|_{\WP}\|x-y\|^{1-\frac{n}{p}},~~~~~x,y\in\Omega.
\ee
Of course the conditions \rf{H-LOC} and \rf{HTOM} with $\|
f\|_{\WPO}$ and  $\|F\|_{\WP}$ replaced by unspecified
constants are not equivalent to membership of $f$ in $\WPO$
or in $\WP|_{\Omega}$ respectively. However the preceding
remarks suggest that a reasonable property which might
perhaps be necessary or perhaps sufficient for a domain
$\Omega $ to be a Sobolev extension domain could be this:
Whenever a function $f:\Omega\rightarrow {\bf R}$ satisfies
$$ |f(x)-f(y)|\le d_{\alpha ,\Omega}
(x,y)^{1-\frac{1}{p}}+\|x-y\|^{1-\frac{n}{p}} $$
for all $x,y\in \Omega $ and $\alpha =(p-n)/(p-1)$ then it
also satisfies
$$ |f(x)-f(y)|\le C(n,p)\|x-y\|^{1-\frac{n}{p}}~ $$
for all $~x,y\in \Omega $ and for some constant $C(n,p)$
depending only on $n$ and $p$.
\par One would like to have a simpler condition on $\Omega
$ which would be sufficient to imply the above ``reasonable
property''. It is clear that the following property, which
has already been considered and studied by other authors,
namely
$$ d_{\alpha ,\Omega }(x,y)^{1-\frac{1}{p}}
\le C\|x-y\|^{1-\frac{n}{p}}
~~\text{ for all}~~x,y\in \Omega ~~\text{and}
~~\alpha =(p-n)/(p-1)
$$
or, equivalently, $d_{\alpha ,\Omega }(x,y)\le C\| x-y\|
^{\alpha }$ for all $x,y\in \Omega $, is such a condition.
\par These considerations lead us to work with a certain
class of domains, essentially those which were introduced
in \cite{GM}. In our context here, it seems convenient to
use terminology different from that of \cite{GM} and other
papers.
\begin{definition}\lbl{ASHD} {\em For each
$\alpha\in(0,1]$, the domain $\Omega \subset {\bf R}^{n}$
is said to be $\alpha $-{\it subhyperbolic} if there exist
constants $C_{\alpha,\Omega}>0$ and
$\theta_{\alpha,\Omega}>0$ such that
$$
d_{\alpha,\Omega }(x,y)\le \,C_{\alpha,\Omega}\| x-y\|
^{\alpha}
$$
for every $x,y\in
\Omega$ satisfying $\| x-y\| \le \theta_{\alpha,\Omega}$.
\par We denote the class of $\alpha $-subhyperbolic domains
in ${\bf R}^{n}$ by $U_{\alpha }({\bf R}^{n})$.}
\end{definition}
\par In \cite{GM} and also in  \cite{L} these domains are
called ``$Lip_{\alpha }$-extension domains". (This name is
derived from the fact that $\Omega\in U_{\alpha }(\RN)$ iff
all functions which are {\it locally} Lipschitz of order
$\alpha$ on $\Omega$ are Lipschitz of order $\alpha$ on
$\Omega$.) These domains have also been studied in
\cite{BSt2,BSt,BSt3} where they are called
``$\alpha-m$cigar domains", and in \cite{BKos} where they
are termed ``local weak $\alpha$-cigar domains".
\par Now Theorem \reff{MAIN-EXT} can be reformulated as
follows: For each $p>n$ and for each
$\frac{p-n}{p-1}$\,-\,sub\-hy\-perbolic domain $\Omega$ in
$\RN$, there exists a constant $\ps\in(n,p)$ depending only
on $n$, $p$ and $\Omega$, such that $\Omega$ is a Sobolev
$W_{q}^{k}$-extension domain for every $q\ge\ps$. \par In
turn, Theorem \reff{RTWO-EXT} admits the following
reformulation: For each $p>2$, a finitely connected bounded
domain $\Omega\subset\R^2$ is a Sobolev
$W_{p}^{1}$-extension domain if and only if $\Omega$ is a
$\frac{p-2}{p-1}$\,-\,subhyperbolic domain.
\par The family $\{U_{\alpha }({\bf R}^{n}): \alpha \in
(0,1]\}$ is an ``increasing family'', i.e.,
$$ U_{\alpha'}({\bf R}^{n})\subset U_{\alpha'' } ({\bf
R}^{n})~~~\text{whenever}~~0<\alpha' <\alpha'' \le 1, $$
see, e.g. \cite{BKos}. Lappalainen \cite{L} proved that
$$
U_{\alpha }({\bf R}^{n})\subsetneqq \bigcap_{\alpha
<\tau\le 1}U_{\tau}(\RN)\text{ for every }\alpha \in (0,1).
$$
\par This last result motivates our discussion presented in
Section 2, which is devoted to the following question:
Does the
equality
\bel{UP-IN} U_{\alpha }(\RN)=\bigcup_{0<\tau <\alpha
}U_{\tau}(\RN) \ee
hold? In other words, do $\alpha$-subhyperbolic domains
have the following ``self-improving" property that whenever
$\Omega$ is an $\alpha$-subhyperbolic domain in $\RN$ for
some $\alpha \in (0,1)$, it is also $\tau$-subhyperbolic
for some positive $\tau$ which is {\it strictly less} than
$\alpha$? (Of course, $\tau$ can depend on $\Omega$).
\par We do not know the answer to this question in general.
We do know that the answer is affirmative for an arbitrary
finitely connected bounded domain $\Omega\in
U_\alpha(\R^2)$, $\alpha\in(0,1)$, as it follows from
Theorem \reff{MAIN-EXT} and Theorem \reff{RTWO-EXT}. We
also know that for a certain subfamily of $U_\alpha(\RN)$,
the so-called {\it strongly $\alpha$-subhyperbolic} domains
(Definition \reff{ST-SD}) the answer to the above question
is affirmative. (See Proposition \reff{CG-D}.) It should be
pointed out that we have no examples of subhyperbolic
domains which are not strongly subhyperbolic.
\par We are able to show that the following weaker version
of the self-improving property \rf{UP-IN} holds for an {\it
arbitrary} subhyperbolic domain in $\RN$.
\begin{theorem}\lbl{S-IM} Let $\alpha\in(0,1)$ and let
$\Omega$ be an $\alpha$-subhyperbolic domain in $\RN$.
There exist a constant $\as,0<\as<\alpha,$ and constants
$\theta,C>0$ such that the following is true:
\par For every $\ve>0$ and every $x,y\in\Omega$,
$\|x-y\|\le\theta,$ there exist a rectifiable curve
$\Gamma\subset\Omega$ joining $x$ to $y$ and a subset
$\GW\subset\Gamma$ consisting of a finite number of arcs
such that the following conditions are satisfied:
\par (i). For every $\tau\in[\as,\alpha]$
\bel{H-GW} \intl_{\GW} \dist(z,\partial\Omega)^{\tau
-1}\,ds(z)\le C\|x-y\|^{\tau}. \ee
In addition, for every ball $B$ centered in $\GW$ of radius
at most $\|x-y\|$,
\bel{H-REG} \diam B\le C\,\lng(B\cap\GW). \ee
\par (ii). We have $\lng(\Gamma)\le C\|x-y\|$ and
\bel{H-LGW} \lng(\Gamma\setminus\GW)<\ve.\ee
Moreover,
\bel{H-AA}
\intl_{\Gamma\setminus\GW}
\dist(z,\partial\Omega) ^{\alpha-1}\,ds(z)\le
C\|x-y\|^{\alpha}.
\ee
\par The constants $\as,\theta$ and $C$ depend only on
$n$, $\alpha$, and the constants $C_{\alpha,\Omega}$ and
$\theta_{\alpha,\Omega}$ introduced in Definition
\reff{ASHD}.
\end{theorem}
\par The proof of this result, presented in Section 2, is
based on the reverse H\"{o}lder inequa\-li\-ty for
$m$-dyadic $A_{1}$-weights. (See Melas \cite{Mel}.)
\par Theorem \reff{S-IM} is an important ingredient in the
proof of the extension Theorem \reff{MAIN-EXT}. It enables
us to prove the following version of the {\it
Sobolev-Poincar\'e inequality} for subhyperbolic domains
(for $p>n$ and $k\ge 1$): Let $\Omega$ be an
$\alpha$-subhyperbolic domain in $\RN$ with
$\alpha=(p-n)/(p-1).$ Given $f\in C^{k-1}(\Omega)$ and
$x\in\Omega$ we let $T_{x}^{k-1}(f)$ denote the Taylor
polynomial of $f$ at $x$ of degree at most $k-1$. We prove
that there exists $\ps\in(n,p)$ and constants
$\theta,\lambda,C>0$ such that for every function $f\in
C^{k-1}(\Omega)\cap\WKPO$ and every $x,y\in\Omega,
\|x-y\|\le\theta,$ the following inequality
$$
|f(y)-T_{x}^{k-1}(f)(y)|\le C\|x-y\|
^{k-\tfrac{n}{\ps}}\left(\,\,\intl_{B\cap\Omega}
 \|\nabla^k f\|^\ps\,dx\right)^{\frac1\ps}
$$
holds. Here $B=B(x,\lambda\|x-y\|)$ is the ball centered at
$x$ of radius $r=\lambda\|x-y\|$. This inequality is a
particular case of Theorem \reff{SH-TE} which we prove in
Section 3.
\par In Section 4 we prove a corollary of this result
related to the sharp maximal function
$$
f^{\sharp}_{k,\Omega}(x):=\sup_{r>0}
r^{-k}\,\inf_{P\in\Pc_{k-1}}
\frac{1}{|B(x,r)|} \intl_{B(x,r)\cap \Omega}
|f-P|\,dx,\ \ \ \ \ \ x\in \Omega. $$
Here $\Pc_{k-1}$ is the space of polynomials of degree at
most $k-1$ defined on $\RN$ and $|B(x,r)|$ is the Lebesgue
measure of the ball $B(x,r)$. We show that for every
$f\in\WKPO$ and every $x\in\Omega$ the following inequality
\bel{E-SF} f^{\sharp}_{k,\Omega}(x)\le
C\left\{(\Mc[(\|\nabla^kf\|^\cw)^{p^*}](x))
^{\frac{1}{p^*}}+\Mc[f^\cw](x)\right\} \ee
holds. Here $\Mc$ denotes the Hardy-Littlewood maximal
function and the symbol $g^\cw$ stands for the extension
{\it by zero} of a function from $\Omega$ to all of $\RN$.
\par The sharp maximal function is a useful tool in the
study of Sobolev functions. In \cite{C} Calder\'{o}n proved
that, for $p>1$, a function $f$ is in $\WKP$ if and only if
$f$ and $f^{\sharp}_{k,\RN}$ are both in $L_{p}(\RN)$. In
\cite{S4} this description has been generalized to the case
of the so-called {\it regular} subsets of $\RN$, i.e., the
sets $S$ such that $|B\cap S|\sim |B|$ for all balls $B$
centered in $S$ of radius at most $1$. We proved in
\cite{S4} that if $S$ is regular and $f\in L_p(S),p>1,$
then $f$ can be extended to a function $F\in \WKP$ if and
only if its sharp maximal function $f_{k,S}^{\sharp}\in
L_p(S).$ (For the case $k=1$ see also \cite{S3,HKT,HKT1}.
Observe that every Sobolev $W^1_p$-extension domain is a
regular subset of $\RN$, see Hajlasz, Koskela and Tuominen
\cite{HKT}. In \cite{S7} we present a description of the
trace space $\WP|_S,p>n,$ for an {\it arbitrary} set
$S\subset\RN$ via an $L_\infty$-version of the sharp
maximal function).
\par Every subhyperbolic domain is a regular set, as shown
in Lemma \reff{BP}. So, in order to prove, for some given
$q>\ps$, that a function $f\in \WKQO$ extends to a Sobolev
$W^k_q$-function on $\RN$, it suffices to show that
$f^{\sharp}_{k,\Omega}\in L_q(\Omega)$. We do this by
applying the Hardy-Littlewood maximal theorem to inequality
\rf{E-SF}. This gives us the inequality
$\|f^{\sharp}_{k,\Omega}\|_{L_q(\Omega)}\le C\|f\|_{\WKQO}$
which completes the proof of Theorem \reff{MAIN-EXT}.
\par {\bf Acknowledgement.} The author is greatly indebted
to Michael Cwikel, Charles Fefferman, Vladimir Maz'ya  and
Naum Zobin for interesting discussions and helpful
suggestions and remarks.
\SECT{2. Subhyperbolic domains: intrinsic metrics and
self-improvement.}{2}
\indent\par Throughout the paper $C,C_1,C_2,...$ will be
generic positive constants which depend only on parameters
determining sets (say, $n,\alpha,$ the constants
$C_{\alpha,\Omega}$ or $\theta_{\alpha,\Omega}$, etc.) or
function spaces ($p,q,$ etc). These constants can change
even in a single string of estimates. The dependence of a
constant on certain parameters is expressed, for example,
by the notation $C=C(n,p)$.
\par The Lebesgue measure of a measurable set $A\subset
\RN$ will be denoted by $\left|A\right|$. Given subsets
$A,B\subset \RN$,  we put
$ \diam A:=\sup\{\|a-a'\|:~a,a'\in A\} $
and
$$
\dist(A,B):=\inf\{\|a-b\|:~a\in A, b\in B\}.
$$
For $x\in \RN$ we also set $\dist(x,A):=\dist(\{x\},A)$.
\par Let $\gamma:[a,b]\to\RN$ be a curve in $\RN$, and let
$u=\gamma(t_1),v=\gamma(t_2)$ where $a\le t_1<t_2\le b.$ By
$\gamma_{uv}$ we denote the arc of $\gamma$ joining $u$ to
$v$.
\par We will be needed the following auxiliary lemma.
\begin{lemma}\lbl{C-G1} (i). Let $x,y\in \Omega$ and let
\bel{ML1}
\max(\dist(x,\partial\Omega),\dist(y,\partial\Omega))\le
2\|x-y\|. \ee
Let $\gamma$ be a rectifiable curve joining $x$ to $y$ in
$\Omega$. Assume that for some $\alpha\in(0,1)$ and $C>0$
the following inequality
\bel{AB-N1} \intl_\gamma\dist(z,\partial
\Omega)^{\alpha-1}\,ds(z) \le C\|x-y\|^\alpha \ee
holds. Then $$\lng(\gamma)\le 2e^C\|x-y\|.$$
\par (ii). Let $x,y\in \Omega$ and let
\bel{ML2}
\max(\dist(x,\partial\Omega),\dist(y,\partial\Omega))>
2\|x-y\|. \ee
Then the line segment $[x,y]\subset \Omega$ and for every
$\beta\in(0,1]$ we have
\bel{LS} \intl_{[x,y]}\dist(z,\partial
\Omega)^{\beta-1}\,ds(z) \le \|x-y\|^\beta. \ee
\end{lemma}
\par {\it Proof.} (i). Let us parameterize $\gamma$ by
arclength; thus we identify $\gamma$ with a function
$\gamma:[0,\ell]\to\Omega$ satisfying $\gamma(0)=x,
\gamma(\ell)=y$. Now \rf{AB-N1} is equivalent to
$$
\intl_0^{\ell}\dist(\gamma(t),\partial
\Omega)^{\alpha-1}\,dt \le C\|x-y\|^\alpha.
$$
Since $\dist(\cdot,\partial\Omega)$ is a Lipschitz function
on $\RN$,
\bel{D-LIP}
\dist(u,\partial\Omega)\le
\dist(v,\partial\Omega)+\|u-v\|,~~~~u,v\in\Omega,
\ee
so that for every $t\in(0,\ell]$
$$ \dist(\gamma(t),\partial\Omega)\le
\dist(x,\partial\Omega)+\|x-\gamma(t)\|. $$
Since $\gamma$ is parameterized by arclength,
$$ \|x-\gamma(t)\|\le \lng(\gamma_{x\gamma(t)})=t, $$
so that
$$
\dist(\gamma(t),\partial\Omega)\le
\dist(x,\partial\Omega)+t,~~~~t\in[0,\ell].
$$
This inequality and \rf{AB-N1} imply
\be
C\|x-y\|^\alpha&\ge&\intl_0^{\ell}\dist(\gamma(t),\partial
\Omega)^{\alpha-1}\,dt \ge\intl_0^{\ell}
(\dist(x,\partial\Omega)+t)^{\alpha-1}\,dt\nn\\&=&
\alpha^{-1}((\dist(x,\partial\Omega)+\ell)^{\alpha}
-\dist(x,\partial\Omega)^{\alpha})\nn\\&\ge&
\alpha^{-1}(\ell^{\alpha}
-\dist(x,\partial\Omega)^{\alpha}).
\nn\ee
But $2\|x-y\|>\dist(x,\partial\Omega)$ so that
$$
C\|x-y\|^\alpha\ge
\alpha^{-1}(\ell^{\alpha}-(2\|x-y\|)^{\alpha}).
$$
Hence
$$
\ell\le (\alpha C+2^\alpha)^{\frac1\alpha}\|x-y\|
\le 2e^C\|x-y\|
$$
proving (i).
\par (ii). Clearly, \rf{ML2} implies $[x,y]\subset \Omega$.
Prove \rf{LS}.
\par We may assume that $\dist(x,\partial\Omega)>
2\|x-y\|.$ Also note that $\|x-z\|\le \|x-y\|$ for every
$z\in[x,y]$. These inequalities and \rf{D-LIP} imply the
following:
$$
\tfrac12\dist(x,\partial\Omega)
\le \dist(x,\partial\Omega)-\|x-y\|
\le \dist(x,\partial\Omega)-\|x-z\|
\le \dist(z,\partial\Omega).
$$
Hence,
\be \intl_{[x,y]}\dist(z,\partial \Omega)^{\beta-1}\,ds(z)
&\le& \intl_{[x,y]}2^{1-\beta}\dist(x,\partial
\Omega)^{\beta-1}\,ds(z)\nn\\
&=&2^{1-\beta}\|x-y\|\dist(x,\partial
\Omega)^{\beta-1}\nn\\
&\le& 2^{1-\beta}\|x-y\| (2\|x-y\|)^{\beta-1}=
\|x-y\|^{\beta}\nn \ee
proving the lemma.\bx
\begin{lemma}\lbl{C-G2} Let $x,y\in \Omega$ and let
$\gamma\subset \Omega$ be a rectifiable curve joining $x$
to $y$. Suppose that for some $\alpha\in(0,1)$ and $C\ge 1$
the following inequality
\bel{AB-N}
\intl_\gamma\dist(z,\partial
\Omega)^{\alpha-1}\,ds(z) \le C\lng^\alpha(\gamma)
\ee
holds. Then
\par (i). There exists a point $\bz\in\gamma$ such that
$$
\lng(\gamma)
\le C^{\frac{1} {1-\alpha}}\,\dist(\bz,\partial\Omega)
$$
\par (ii).  We have
$$
\frac{1}{\lng(\gamma)}\intl_\gamma \dist(z,\partial
\Omega)^{\alpha-1}\,ds(z)\le 2C\,\inf_{z\in\gamma}
\dist(z,\partial \Omega)^{\alpha-1}.
$$
\end{lemma}
\par {\it Proof.} (i). Put $\ell:=\lng(\gamma)$. Let $\bz$
be a point in $\gamma$ such that
$$
\max\{\dist(z,\partial\Omega):~z\in\gamma\}=
\dist(\bz,\partial\Omega).
$$
Then
$$
\intl_\gamma\dist(z,\partial
\Omega)^{\alpha-1}\,ds(z)
\ge\intl_\gamma
\dist(\bz,\partial\Omega)^{\alpha-1}\,ds(z)=
\ell\dist(\bz,\partial\Omega)^{\alpha-1}
$$
so that, by \rf{AB-N},
$$
\ell\dist(z,\partial\Omega)^{\alpha-1}\le
C\ell^\alpha.
$$
Hence
$$
\ell\le C^{\frac{1} {1-\alpha}}\dist(z,\partial\Omega)
$$
proving (i).
\par (ii). Put $ w(z):=\dist(z,\partial\Omega).$ Then, by
\rf{AB-N},
\bel{H-PN}
\frac{1}{\ell}\intl_\gamma
w(z)^{\alpha-1}\,ds(z) \le \ell^{-1}(C\,\,\ell^\alpha)=
C\,\ell^{\alpha-1}.
\ee
For every $z_1,z_2\in\gamma$ we have
$$
|w(z_1)-w(z_2)|=
|\dist(z_1,\partial\Omega)-\dist(z_2,\partial\Omega)|\le
\|z_1-z_2\|\le \ell
$$
so that
\bel{L-DMN}
\max_{z\in\gamma}w(z)\le
\min_{z\in\gamma}w(z)+\ell.
\ee
\par Let us consider two cases. First suppose that
$\max_{z\in\gamma} w(z)\le 2\ell.$ Since $\alpha\in(0,1)$,
we obtain
$$
\ell^{\alpha-1}\le 2^{1-\alpha}
\min_{z\in\gamma}w(z)^{\alpha-1}
$$
so that, by \rf{H-PN},
$$
\frac{1}{\ell}\intl_\gamma
w(z)^{\alpha-1}\,ds(z)
\le 2^{1-\alpha}C \min_{z\in\gamma}w(z)^{\alpha-1}\le
2C \min_{z\in\gamma}w(z)^{\alpha-1}.
$$
\par Now assume that $2\ell<\max_{z\in\gamma} w(z).$ Then,
by \rf{L-DMN},
$$
\max_{z\in\gamma} w(z)\le \min_{z\in\gamma}w(z)+
\tfrac12\max_{z\in\gamma} w(z)
$$
so that $\max_{z\in\gamma}w(z)\le 2\min_{z\in\gamma}w(z).$
Hence
$$
\max_{z\in\gamma} w(z)^{\alpha-1}\le 2^{1-\alpha}
\min_{z\in\gamma}w(z)^{\alpha-1}\le 2
\min_{z\in\gamma}w(z)^{\alpha-1}.
$$
Finally, we have
$$
\frac{1}{\ell}\intl_\gamma w(z)^{\alpha-1}\,ds(z)
\le\max_{z\in\gamma} w(z)^{\alpha-1}\le
2\min_{z\in\gamma}w(z)^{\alpha-1}.
$$
\par The lemma is proved.
\par This lemma implies the following important property of
subhyperbolic domains.
\begin{lemma}\lbl{BP} Let $\alpha\in(0,1)$ and let $\Omega$
be an $\alpha$-subhyperbolic domain.
\par There exist constant $\delta>0$ and $\sigma\in(0,1]$
depending only on $n,\alpha, C_{\alpha,\Omega}$ and
$\theta_{\alpha,\Omega}$ such that every ball $B$ centered
in $\Omega$ of diameter at most $\delta$ contains a ball
$B'\subset\Omega$ of diameter at least $\sigma\diam B$.
\end{lemma}
\par {\it Proof.} Let
$\delta:=\min\{\theta_{\alpha,\Omega},\tfrac12\diam
\Omega\}$. Let $B=B(x,r)$ be a ball with center in
$x\in\Omega$ and radius $r\le\delta$. Put
$\tr:=r/(8e^{C_{\alpha,\Omega}})$. Since
$\tr\le\delta\le\tfrac12\diam \Omega,$ there exists a point
$a\in\Omega$ such that $\|x-a\|>\tr$. Let
$\Gamma\subset\Omega$ be a curve joining $x$ to $a$. Since
$a\notin B(x,r)$, we have
$\Gamma\cap\partial(B(x,r))\ne\emp$ so that there exists a
point $b\in\Omega$ such that $\|x-b\|=\tr.$
\par If
$$ \max(\dist(x,\partial\Omega),\dist(b,\partial\Omega))>
2\|x-b\|=2\tr, $$
then either $B(x,\tr)\subset B(x,r)\cap\Omega$ or
$B(b,\tr)\subset B(x,r)\cap\Omega$, so that in this case
the ball $B'$ exists.
\par Suppose that
$$ \max(\dist(x,\partial\Omega),\dist(b,\partial\Omega))\le
2\|x-b\|=2\tr, $$
Since $\tr\le\delta\le\theta_{\alpha,\Omega}$, there exists
a curve $\gamma\subset\Omega$ joining $x$ to $b$ such that
$$
\intl_\gamma\dist(z,\partial \Omega)^{\alpha-1}\,ds(z)
\le C_{\alpha,\Omega}\|x-b\|^\alpha.
$$
(We may assume that $C_{\alpha,\Omega}\ge 1$.) By Lemma
\reff{C-G1}, part (i),
$$
\lng(\gamma)\le 2e^{C_{\alpha,\Omega}}\|x-b\|
=2e^{C_{\alpha,\Omega}}\tr=r/4.
$$
Moreover, by part (ii) of Lemma \reff{C-G2}, there exists a
point $\bz\in\gamma$ such that
$$
\lng(\gamma)\le C_{\alpha,\Omega}^{\frac{1}{1-\alpha}}
\dist(\bz,\partial \Omega).
$$
Hence,
$$
\tr=\|x-b\|\le \lng(\gamma)
\le C_{\alpha,\Omega}^{\frac{1}{1-\alpha}}
\dist(\bz,\partial \Omega).
$$
Put $r':= \tr/(2C_{\alpha,\Omega}^{\frac{1}{1-\alpha}})$
and $B':=B(z,r')$. Then, $r'\le \tr\le r/4$ (recall that
$C_{\alpha,\Omega}\ge 1$) and
$ r'<\dist(\bz,\partial \Omega).$
\par Hence, $B(z,\bar{r})\subset\Omega.$ On the other hand,
$ \|x-z\|\le \lng(\gamma)\le r/4 $
so that $B'\subset B(x,r/2)\subset B$. The lemma is
proved.\bx
\par Before to present the proof of Theorem \reff{S-IM} let
us demonstrate its main ideas for a family of the so-called
strongly subhyperbolic domains in $\RN$
\begin{definition}\lbl{ST-SD} {\em Let $\alpha\in(0,1]$. A
domain $\Omega \subset {\bf R}^{n}$ is said to be {\it
strongly $\alpha $-subhyperbolic} if there exist constants
$C\,,\theta>0$ such that every $x,y\in\Omega$, $\|x-y\|\le
\theta$, can be joined by a rectifiable curve
$\gamma\subset\Omega$ satisfying the following condition:
for every $u,v\in\gamma$
\bel{SG-IH} \int_{\gamma_{uv}}\dist(z,\partial\Omega)
^{\alpha-1}\,ds(z) \le\, C\|u-v\|^{\alpha}. \ee
}\end{definition}
\begin{remark} {\em Given $x,y\in\Omega$ a rectifiable
curve $\gamma\subset\Omega$ joining $x$ to $y$ is said to
be {\it $d_{\alpha,\Omega}$-geodesic} if
$$
d_{\alpha,\Omega}(x,y)=\len_{\alpha ,\Omega }(\gamma
):=
\intl_{\gamma}\dist(z,\partial\Omega)^{\alpha-1}\,ds(z).
$$
(See definition \rf{DEF-D}.)
\par Clearly, if $\Omega$ is $\alpha$-subhyperbolic and for
every $x,y\in\Omega$ there exists
$d_{\alpha,\Omega}$-geodesic, then $\Omega$ is strongly
$\alpha$-subhyperbolic. In fact, in this case every arc of
$d_{\alpha,\Omega}$-geodesic curve is
$d_{\alpha,\Omega}$-geodesic as well so that inequality
\rf{SG-IH} holds.
\par However, for every $\alpha\in(0,1]$ there exists a
domain $\Omega\in\RN$ and $x,y\in\Omega$ such that
$d_{\alpha,\Omega}$-geodesic for $x,y$ does not exist. This
is trivial for $\alpha=1$, i.e., for quasi-Euclidean
domains. For the case $\alpha\in(0,1)$ see \cite{BSt2}.
\par Let us slightly generalize this example. Fix $C\ge 1$.
We say that a rectifiable curve $\gamma\subset\Omega$
joining $x$ to $y$ is {\it
$(C,d_{\alpha,\Omega})$-geodesic} if for every
$u,v\in\gamma$ the following inequality
$$
\len_{\alpha,\Omega }(\gamma_{uv})\le
C\,d_{\alpha,\Omega}(u,v).
$$
holds. Clearly, a rectifiable curve $\gamma$ is
$(1,d_{\alpha,\Omega})$-geodesic iff it is
$d_{\alpha,\Omega}$-geodesic. Moreover, if $\Omega\in
U_\alpha(\RN)$ and for every $x,y\in \Omega$,
$\|x-y\|\le\theta,$ there exists
$(C,d_{\alpha,\Omega})$-geodesic joining $x$ to $y$ in
$\Omega$, then $\Omega$ is strongly $\alpha$-subhyperbolic.
\par This observation motivates the following question: Let
$\Omega$ be a domain in $\RN$ and let $\alpha\in (0,1]$.
Does there exist a constant $C=C_\Omega>1$ such that every
two points $x,y\in\Omega$ can be joined by a
$(C,d_{\alpha,\Omega})$-geodesic curve? Even for the
quasi-Euclidean domains, i.e., for $\alpha=1$, we do not
know the answer to this question.\bx
}\end{remark}
\begin{proposition}\lbl{CG-D} Let $\alpha\in(0,1)$ and let
$\Omega$ be a strongly $\alpha$-subhyperbolic domain in
$\RN$. Then $\Omega$ is $\tau$-subhyperbolic for some
$\tau\in(0,\alpha)$.
\end{proposition}
\par {\it Proof.} Since $\Omega$ is strongly
$\alpha$-subhyperbolic, there exist constants $\theta>0$
and $C\ge 1$ such that every $x,y\in\Omega$, $\|x-y\|\le
\theta$, can be joined by a rectifiable curve
$\gamma\subset\Omega$ satisfying the following condition:
for every $u,v\in\gamma$
$$
\int_{\gamma_{uv}}\dist(z,\partial\Omega)
^{\alpha-1}\,ds(z) \le\, C\|u-v\|^{\alpha}.
$$
In particular,
\bel{U-XY} \intl_{\gamma}
\dist(z,\partial\Omega)^{\alpha-1}\,ds(z)\le C\,
\|x-y\|^\alpha .\ee
\par Let $\ell:=\lng(\gamma)$. We parameterize $\gamma$ by
arclength: thus $\gamma:[0,\ell]\to\Omega$, $\gamma(0)=x,
\gamma(\ell)=y$.  Let $u=\gamma(t_1),v=\gamma(t_2)$ where
$0<t_1<t_2\le \ell$ (recall that by $\gamma_{uv}$ we denote
the arc of $\gamma$ joining $u$ to $v$).
\par Applying part (ii) of Lemma \reff{C-G2} to the arc
$\gamma_{uv}$, we obtain
$$
\frac{1}{\lng(\gamma_{uv})}\intl_{\gamma_{uv}}
\dist(z,\partial \Omega)^{\alpha-1}\,ds(z)\le
2C\,\inf_{z\in\gamma_{uv}} \dist(z,\partial
\Omega)^{\alpha-1},
$$
or, in the parametric form,
\bel{R-F} \frac{1}{t_2-t_1}\intl_{t_1}^{t_2}
\dist(\gamma(t),\partial \Omega)^{\alpha-1}\,dt\le
2C\,\inf_{t\in[t_1,t_2]}\dist(\gamma(t),\partial
\Omega)^{\alpha-1}. \ee
\par Put
$$
w(t):=\dist(\gamma(t),\partial\Omega).
$$
By \rf{R-F}, the function $w=w(t)$ has the following
property: for every subinterval $I\subset [0,\ell]$
$$ \frac{1}{|I|}\intl_{I} w(t)^{\alpha-1}\,dt\le
2C\,\inf_{I} w(t)^{\alpha-1}. $$
\par Thus the function $h:=w^{\alpha-1}$ is a Muckenhoupt's
$\Ac_1-$ weight on $[0,\ell]$, see , e.g., \cite{GR}, so
that $h$ satisfies the reverse H\"older inequality on
$[0,\ell]$ (see Mackenhoupt \cite{Mac}, Gehring \cite{G},
Coiffman and Fefferman \cite{CF}): There exist constants
$\tq>1$ and $C_1\ge 1$ (depending only on $C$) such that
$$
\left(\frac{1}{\ell}\intl_{0}^\ell
h^{\tq}(t)\,dt\right)^{1/\tq} \le
C_1\,\frac{1}{\ell}\intl_{0}^\ell h(t)\,dt.
$$
\par By \rf{U-XY},
$$
\intl_{0}^\ell w^{\alpha-1}(t)\,dt
\le C\, \|x-y\|^\alpha
$$
so that
$$
\left(\frac{1}{\ell}\intl_{0}^\ell
w^{(\alpha-1)\tq}(t)\,dt\right)^{1/\tq} \le
C_1\,\frac{1}{\ell}\intl_{0}^\ell w^{\alpha-1}(t)\,dt
\le C_1\,C\frac{1}{\ell}\|x-y\|^\alpha.
$$
We put $q:=\min\{\tq,\tfrac{1-\alpha/2}{1-\alpha}\}$ and
$\tau=q(\alpha-1)+1$. Clearly, $1<q\le \tq$ and
$0<\tau<\alpha$. Hence,
\be \left(\frac{1}{\ell}\intl_{0}^\ell
w^{\tau-1}(t)\,dt\right)^{1/q}&=&
\left(\frac{1}{\ell}\intl_{0}^\ell
w^{q(\alpha-1)}(t)\,dt\right)^{1/q}\nn\\&\le&
\left(\frac{1}{\ell}\intl_{0}^\ell
w^{\tq(\alpha-1)}(t)\,dt\right)^{1/\tq}\le
C_2\frac{1}{\ell}\|x-y\|^\alpha\nn \ee
where $C_2:=C_1\,C$. Finally, we obtain
\be \intl_{\gamma}\dist(z,\partial\Omega)^{\tau-1}\,ds(z)
&=&\intl_{0}^\ell w^{\tau-1}(t)\,dt\le
\frac{C_2^{q}}{\ell^{q-1}}\|x-y\|^{\alpha q} \nn\\&\le&
C_2^{q}\,\|x-y\|^{\alpha q -q+1}=C_2^{q}\,\|x-y\|^{\tau}
\nn
\ee
proving that $\Omega$ is a $\tau-$subhyperbolic domain.\bx
\medskip
\par {\it Proof of Theorem \reff{S-IM}.} Let $\ve>0$ and
let $\Omega\in U_\alpha(\RN)$. We will assume that the
constant $C_{\alpha,\Omega}\ge 1$. Put
$$
\theta:=\tfrac12\,e^{-2C_{\alpha,\Omega}}
\theta_{\alpha,\Omega}
$$
and fix $x,y\in\Omega$ such that $\|x-y\|\le \theta$.
\par By part (ii) of Lemma \reff{C-G1}, if inequality
\rf{ML2} is satisfied, then the statement of Theorem
\reff{S-IM} is true with $\Gamma=\GW=[x,y]$ and any
$\as\in(0,\alpha)$.
\par Now suppose that $x,y$ satisfy inequality \rf{ML1},
i.e.,
$$
\max(\dist(x,\partial\Omega),\dist(y,\partial\Omega))\le
2\|x-y\|.
$$
\begin{lemma}\lbl{DL-S} Let $\Omega\in U_{\alpha}(\RN)$ and
let $x,y\in\Omega$, $\|x-y\|\le\theta$. Let $0<\delta\le
d_{\alpha,\Omega}(x,y)$ and let $\gamma\subset\Omega$ be a
rectifiable curve joining $x$ to $y$ such that
\bel{O-G} \int_{\gamma}\dist(z,\partial\Omega)
^{\alpha-1}\,ds(z)< d_{\alpha,\Omega}(x,y)+\delta. \ee
Then:
\par (i). We have
$$ \int_{\gamma}\dist(z,\partial\Omega)
^{\alpha-1}\,ds(z)\le 2C_{\alpha,\Omega}\|x-y\|^\alpha$$
and $ \lng (\gamma)\le 2e^{2C_{\alpha,\Omega}}\|x-y\|.$
\par (ii). For every $u,v\in\gamma$ such that
\bel{L-DL} \lng(\gamma_{uv})\ge \delta^{\frac{1}{\alpha}}
\ee
the following inequality
$$
\int_{\gamma_{uv}}\dist(z,\partial\Omega)
^{\alpha-1}\,ds(z)\le
2C_{\alpha,\Omega}\lng^\alpha(\gamma_{uv})
$$
holds.
\end{lemma}
\par {\it Proof.} First prove that
\bel{OPT-D}
\int_{\gamma_{uv}}\dist(z,\partial\Omega)
^{\alpha-1}\,ds(z)<d_{\alpha,\Omega}(u,v)+\delta.
\ee
In fact, assume that
$$
d_{\alpha,\Omega}(u,v)+\delta\le\int_{\gamma_{uv}}
\dist(z,\partial\Omega)^{\alpha-1}\,ds(z).
$$
Then
\be \int_{\gamma}
\dist(z,\partial\Omega)^{\alpha-1}\,ds(z)
&=&\int_{\gamma_{xu}}
\dist(z,\partial\Omega)^{\alpha-1}\,ds(z)
+\int_{\gamma_{uv}}
\dist(z,\partial\Omega)^{\alpha-1}\,ds(z)\nn\\
&+&\int_{\gamma_{vy}}
\dist(z,\partial\Omega)^{\alpha-1}\,ds(z)\nn\\
&\ge&
d_{\alpha,\Omega}(x,u)+
(d_{\alpha,\Omega}(u,v)+\delta)+d_{\alpha,\Omega}(v,y)
\nn
\ee
so that, by the triangle inequality for the metric
$d_{\alpha,\Omega}$,
$$
\int_{\gamma}
\dist(z,\partial\Omega)^{\alpha-1}\,ds(z)\ge
d_{\alpha,\Omega}(x,y)+\delta
$$
which contradicts inequality \rf{O-G}.
\par Since $0<\delta\le d_{\alpha,\Omega}(x,y)$, by
\rf{O-G},
$$ \int_{\gamma}\dist(z,\partial\Omega) ^{\alpha-1}\,ds(z)<
2d_{\alpha,\Omega}(x,y).$$
Since $\theta\le \theta_{\alpha,\Omega}$ and $\Omega\in
U_\alpha(\RN)$, we have $d_{\alpha,\Omega}(x,y)\le
C_{\alpha,\Omega}\|x-y\|^\alpha$ so that
$$
\int_{\gamma}\dist(z,\partial\Omega) ^{\alpha-1}\,ds(z)<
2C_{\alpha,\Omega}\|x-y\|^\alpha.
$$
By Lemma \reff{C-G1}, part (i),
$$
\lng (\gamma)\le 2e^{2C_{\alpha,\Omega}}\|x-y\|
$$
proving (i).
Hence,
$$
\lng (\gamma)\le
2e^{2C_{\alpha,\Omega}}\theta=2e^{2C_{\alpha,\Omega}}
(\tfrac12e^{-2C_{\alpha,\Omega}})\theta_{\alpha,\Omega}=
\theta_{\alpha,\Omega} $$
so that for every $u,v\in\gamma$ we have $\|u-v\|\le\lng(
\gamma)\le\theta_{\alpha,\Omega}$. Since $\Omega\in
U_\alpha(\RN)$, this implies
$$
d_{\alpha,\Omega}(u,v)\le\, C_{\alpha,\Omega}\|u-v\|^\alpha
\le \, C_{\alpha,\Omega}\lng^\alpha(\gamma_{uv}).
$$
Combining this inequality with \rf{OPT-D} and \rf{L-DL}, we
obtain
\be \int_{\gamma_{uv}}\dist(z,\partial\Omega)
^{\alpha-1}\,ds(z)&\le&
C_{\alpha,\Omega}\lng^\alpha(\gamma_{uv}) +\delta\nn\\
&\le& C_{\alpha,\Omega}\lng^\alpha(\gamma_{uv})
+\lng^\alpha(\gamma_{uv})\nn\ee
proving (ii) and the lemma.\bx
\par Put
\bel{D-M}
m:=[2(2C_{\alpha,\Omega})^{\frac{1}{1-\alpha}}]+1. \ee
\par Let $k$ be a positive integer such that
\bel{DEF-K}
2e^{2C_{\alpha,\Omega}}\|x-y\|(1-1/m)^k\le
\ve.
\ee
Finally, we put
$$
\delta:=\min\{d_{\alpha,\Omega}(x,y),
m^{-\alpha k}\|x-y\|^\alpha\}.
$$
Thus $0<\delta\le d_{\alpha,\Omega}(x,y)$ and
\bel{DEF-DL} \|x-y\|m^{-k}\ge\delta^{\frac{1}{\alpha}}.\ee
\par Let $\Gamma\subset\Omega$ be a rectifiable curve
joining $x$ to $y$ such that
$$
\int_{\Gamma}\dist(z,\partial\Omega) ^{\alpha-1}\,ds(z)<
d_{\alpha,\Omega}(x,y)+\delta.
$$
Then, by Lemma \reff{DL-S},
\bel{PR-C}
\int_{\Gamma}\dist(z,\partial\Omega)
^{\alpha-1}\,ds(z)\le 2C_{\alpha,\Omega}\|x-y\|^\alpha
\ee
and
\bel{LGM}
\lng (\Gamma)\le 2e^{2C_{\alpha,\Omega}}\|x-y\|.
\ee
Moreover, for every $u,v\in\Gamma$ such that
\bel{L-DL1} \lng(\Gamma_{uv})\ge\delta^{\frac{1}{\alpha}}
\ee
the following inequality
\bel{D-31} \int_{\Gamma_{uv}}\dist(z,\partial\Omega)
^{\alpha-1}\,ds(z)\le
2C_{\alpha,\Omega}\lng^\alpha(\Gamma_{uv})
\ee
holds.
\par Inequality \rf{D-31} and part (ii) of Lemma
\reff{C-G2} imply:
\bel{R-EH} \frac{1}{\lng (\Gamma_{uv})}\intl_{\Gamma_{uv}}
\dist(z,\partial \Omega)^{\alpha-1}\,ds(z)\le
4C_{\alpha,\Omega}\,\inf_{z\in\Gamma_{uv}} \dist(z,\partial
\Omega)^{\alpha-1}. \ee
\par Put
$$
L:=\lng (\Gamma).
$$
Since $\|x-y\|\le L$, by \rf{DEF-DL} and \rf{L-DL1}, for
every $u,v\in\Gamma$ such that
\bel{L-DL2} \lng (\Gamma_{uv})\ge L\,m^{-k} \ee
inequality \rf{D-31} is satisfied.
\par By $\IM$ we denote the family of all $m$-adic closed
subintervals of the interval $I_0:=[0,L]$. recall that this
family of intervals can be obtained by the standard
iterative procedure: we start with the entire interval
$[0,L]$ and, at each level of the construction, we split
every interval of the given level into $m$ equally sized
closed subintervals.
\par Let
$$
\Ic_{j,m}:=\{{\rm all}~
m{\rm-adic~intervals~of~the~} j-{\rm th~level}\}.
$$
Thus $\Ic_{0,m}:=\{[0,L]\}$,
$\Ic_{1,m}:=\{[Li/m,L(i+1)/m]:~i=0,1,...,m-1\}$ etc.
Clearly, $|I|=L\,m^{-j}$ for every interval $I$ of the
$j$-th level. Put
$$
S_{k,m}:=\bigcup_{j=0}^k\Ic_{j,m}=\{{\rm all}~
m{\rm-adic~intervals~of~the~level~at~most}~~k\}.
$$
\par Let us parameterize $\Gamma$ by arclength; thus we
identify $\gamma$ with a function $\Gamma:[0,L]\to\Omega$
satisfying $\Gamma(0)=x, \Gamma(L)=y$. Finally, put
$$
g(t):= \dist(\Gamma(t),\partial
\Omega)^{\alpha-1},~~~t\in[0,L].
$$
Then, by \rf{L-DL2} and \rf{R-EH}, the following is true:
For each $m$-adic interval $I\in S_{k,m}$ we have
\bel{A1-W} \frac{1}{|I|}\intl_{I} g(t)\,dt\le
C_g\,\inf_{t\in I} g(t) \ee
with $C_g:=4C_{\alpha,\Omega}$.
\par Following Melas \cite{Mel} we say that $g$ is a
Muckenhoupt $A_1$-weight on $[0,L]$ with respect to the
family $S:=S_{k,m}$ of all $m$-adic intervals of the level
at most $k$. We let $\Mc_S$ denote the corresponding
maximal operator for the family $S$:
\bel{DMM} \Mc_S\,g(t):=\sup\left\{\frac{1}{|I|}\intl_{I}
|g(t)|\,dt: I\ni t, I\in S\right\}. \ee
Thus \rf{A1-W} is equivalent to the inequality
$$
\Mc_S\,g(t)\le C_g\, g(t),~~~t\in[0,L].
$$
\par Put
$$
q^\sharp:=\frac{\log m}{\log\left(m-(m-1)/C_g\right)}
$$
and $q^*:=(1+q^\sharp)/2$. Clearly, $1\le C_g<\infty$ so
that $q^\sharp, q^*>1$.
\par We will be needed the following corollary of a general
result proved in \cite{Mel}.
\begin{theorem} For any $A_1$-weight $g$ (with respect to
$S$) and any $q$, $1\le q\le q^*$, the following inequality
\bel{TM}
\left(\frac{1}{L}\intl_0^L(\Mc_S\,g)^q\,dt\right)^{\frac1q}
\le \tC\left(\frac{1}{L}\intl_0^L g\,dt\right) \ee
holds. Here $\tC$ is a constant depending only on $m$ and
$C_g$.
\end{theorem}
\begin{remark} {\em Actually the theorem is true for $1\le
q< q^\sharp$ but with  $\tC$ depending  on $m,C_g$ and $q$,
see \cite{Mel}.}
\end{remark}
\begin{corollary}\lbl{CORR} For any $A_1$-weight $g$ (with
respect to $S$), any family $\Ac$ of non-overlapping
$m$-adic intervals of the level at most $k$ and any $q$,
$1\le q\le q^*$, we have
\bel{C-ST}
\left(\frac{1}{L}\sum_{I\in\Ac}\left(\frac{1}{|I|} \intl_I
g\,dt\right)^q|I|\right)^{\frac1q} \le
\tC\left(\frac{1}{L}\intl_0^L g\,dt\right)\,.
\ee
\end{corollary}
\par {\it Proof.} In fact, by definition \rf{DMM}, for
every $I\in\Ac$ and every $t\in I$
$$
\frac{1}{|I|} \intl_I g\,ds\le \Mc_S\,g(t)
$$
so that
$$
\left(\frac{1}{|I|} \intl_I
g\,dt\right)^q|I|\le \intl_I(\Mc_S\,g)^q\,dt.
$$
Therefore the left-hand side of \rf{C-ST} does not exceed
$$
\left(\frac{1}{L}\sum_{I\in\Ac}\intl_I(\Mc_S\,g)^q\,dt
\right)^{\frac1q} \le
\left(\frac{1}{L}\intl_0^L(\Mc_S\,g)^q\,dt\right)
^{\frac1q}
$$
which together with \rf{TM} implies the required inequality
\rf{C-ST}.\bx
\par We turn to construction of a family $\Ac\subset
S_{k,m}$ of non-overlapping $m$-adic intervals of the level
at most $k$ such that for each $I\in\Ac$
$$ \sup_I g\le C\inf_I g $$
and
$$
|\,[0,L]\setminus\{\cup I:~I\in\Ac\}|<\ve.
$$
Here $C$ is a constant depending only on $n,\alpha,$ and
$C_{\alpha,\Omega}$.
\par Let $I=[t_1,t_2]\in S_{k-1,m}$ be an $m$-adic interval
of the level at most $k-1$ and let $u:=\Gamma(t_1),
v:=\Gamma(t_2)$.  By \rf{D-31} and part (i) of Lemma
\reff{C-G2}, there exists $t_I\in I$ such the point
$z_I=\Gamma(t_I)\in\Gamma_{uv}$ satisfies the following
inequality:
\bel{UV-Z} \lng(\Gamma_{uv})\le \,
C'\dist(z_I,\partial\Omega) \ee
with $C':=(2C_{\alpha,\Omega})^{\frac{1}{1-\alpha}}.$
\par Let us split the interval $I$ into $m$ equal
subintervals $I^{(1)},...,I^{(m)}$. Then $t_I\in I^{(j)}$
for some $j\in\{1,...,m\}$. By
$\Gamma^{(j)}:=\Gamma|_{I^{(j)}}$ we denote the arc
corresponding to the interval $I^{(j)}$. Thus $I^{(j)}\in
S_{k,m}$ is an $m$-adic interval of the level at most $k$
and
$$
\lng (\Gamma^{(j)})=|I^{(j)}|=|I|/m
=\lng (\Gamma_{uv})/m.
$$
\par Since $\dist(\cdot,\partial\Omega)$ is a Lipschitz
function,  for every $t\in I^{(j)}$ we have
$$
|\dist(z_I,\partial\Omega)-\dist(z(t),\partial\Omega)|\le
\|z_I-z(t)\|\le \lng (\Gamma^{(j)})
=\lng(\Gamma_{uv})/m.
$$
Combining this inequality with \rf{UV-Z} we obtain
$$
|\dist(z_I,\partial\Omega)-\dist(z(t),\partial\Omega)|\le
\frac{C'}{m}\dist(z_I,\partial\Omega).
$$
But $m:=[2C']+1$, see \rf{D-M}, so that
$$
|\dist(z_I,\partial\Omega)-\dist(z(t),\partial\Omega)|\le
\tfrac12\dist(z_I,\partial\Omega).
$$
Hence
$$
\tfrac12\dist(z_I,\partial\Omega)
\le\dist(z(t),\partial\Omega)\le
\tfrac32\dist(z_I,\partial\Omega),~~~t\in I^{(j)}.
$$
\par We let $\tI$ denote the interval $I^{(j)}$. Thus we
have proved that for each $I\in S_{k-1,m}$ there exists a
subinterval $\tI\in S_{k,m}$, $\tI\subset I$, such that
$$
\max_{t\in\tI}\dist(z(t),\partial\Omega)\le
3\min_{t\in\tI}\dist(z(t),\partial\Omega).
$$
Since $g(t):=\dist(z(t),\partial\Omega)^{\alpha-1}$, we
obtain
$$
\max_{t\in\tI} g(t)\le
3^{1-\alpha}\min_{t\in\tI}g(t).
$$
\par Now we construct the family $\Ac\subset S_{k,m}$ as
follows. At the first stage for the interval $I_0:=[0,L]$
we determine an $m$-adic interval $\tI_0\in \Ic_{1,m}$ of
the first level and put $\Ac_1:=\{\tI_0\}$ and
$U_1:=\tI_0$.
\par Let us consider the set $[0,L]\setminus U_1$ which
consists of $m-1$ $m$-adic intervals of the first level. We
let $\Bc_1$ denote the family of these intervals. For every
$I\in\Bc_1$ we construct the interval $\tI\in \Ic_{2,m}$
and put
$$
\Ac_2:=\{\tI\in \Ic_{2,m}:~I\in\Bc_1\}.
$$
By $U_2$ we denote the set
$$
U_2:=U_1\cup\{I:~I\in\Ac_2\}.
$$
\par Now the set $[0,L]\setminus U_2$ consists of $(m-1)^2$
$m$-adic intervals of the second level. We denote the
family of these intervals by $\Bc_2$ and finish the second
stage of the procedure.
\par After the $k$-th stages of this procedure we obtain
the families $\Ac_j\subset\Ic_{j,m}$, $j=1,2,...,k,$ of
$m$-adic intervals. We put
$$
\Ac=\cup\{\Ac_j:~j=1,...,k\}.
$$
Thus $\Ac\subset S_{k,m}$ is a family of $m$-adic intervals
of the level at most $k$. We know that for every interval
$I\in \Ac$ the following inequality
\bel{EQ-U} \max_{t\in I} g(t)\le 3^{1-\alpha}\min_{t\in
I}g(t)\ee
holds. We also know that the set
$$
U=U_k:=\cup\{I:~I\in\Ac\}
$$
has the following property: the set
$$
E:=[0,L]\setminus U=[0,L]\setminus\cup\{I:~I\in\Ac\}
$$
consists of $(m-1)^k$ $m$-adic intervals of the $k$-th
level. Since $|I|=m^{-k}\,L$ for each $I\in\Ic_{k,m}$, we
obtain
$$
|E|=\frac{(m-1)^k}{m^k}\,L.
$$
But, by \rf{LGM},
$$ L=\lng (\Gamma)\le C''\,\|x-y\|$$
where $C'':=2e^{2C_{\alpha,\Omega}}$. Hence,
$$
|E|\le\,C''\,\|x-y\|(1-1/m)^k.
$$
Combining this inequality with \rf{DEF-K}, we obtain the
required estimate
\bel{F-B} |E|=|\,[0,L]\setminus U|\le \ve. \ee
\par Now for the family $\Ac$ constructed above let us
estimate from below the quantity
$$
T:=\left(\frac{1}{L}\sum_{I\in\Ac}\left(\frac{1}{|I|}
\intl_I g\,dt\right)^q|I|\right)^{\frac1q}
$$
from the left-hand side of inequality \rf{C-ST}. By
\rf{EQ-U}, for each $I\in\Ac$ we have
$$
\intl_I g^q\,dt\le |I|\,\max_{t\in I} g^q(t)\le
3^{q(1-\alpha)}|I|\,(\min_{t\in I}g(t))^q
\le 3^{q(1-\alpha)}|I|
\left(\frac{1}{|I|}\intl_I g\,dt\right)^q
$$
so that
$$
T^q\ge 3^{q(\alpha-1)}\,\frac{1}{L}\sum_{I\in\Ac}
\intl_I g^q\,dt=3^{q(\alpha-1)}\,\frac{1}{L}
\intl_U g^q\,dt.
$$
(Recall that $U=\cup\{I:~I\in\Ac\}$.)
\par By Corollary \reff{CORR},
$$ T\le \tC\left(\frac{1}{L}\intl_0^L g\,dt\right) $$
so that
$$
\frac{1}{L}
\intl_U g^q\,dt\le 3^{q(1-\alpha)}T^q\le
C_1\left(\frac{1}{L}\intl_0^L g\,dt\right)^q\,.
$$
with $C_1:=3^{q(\alpha-1)} \tC$. On the other hand, by
inequality \rf{PR-C},
$$
\intl_0^L g\,dt=\int_{\Gamma}\dist(z,\partial\Omega)
^{\alpha-1}\,ds(z)\le 2C_{\alpha,\Omega}\|x-y\|^\alpha.
$$
Hence
$$
\frac{1}{L}
\intl_U g^q\,dt\le
C_1\left(\frac{1}{L}\intl_0^L g\,dt\right)^q\le
C_2\|x-y\|^{q\alpha}/L^q
$$
with $C_2:=(2C_{\alpha,\Omega})^q C_1$.
\par Recall that this inequality holds for all $q\in
[1,q^*]$, see Corollary \reff{CORR}. We put
$\tq:=\min\{q^*,\frac{1-\alpha/2}{1-\alpha}\}$ and
$\as:=1-\tq(1-\alpha)$. Since $q^*>1$ and $0<\alpha<1$, we
have $1<\tq\le q$ and $0<\as<\alpha$.
\par Let $\tau\in[\as,\alpha]$ and let
$q:=\frac{1-\tau}{1-\alpha}$. Then $q\in [1,q^*]$ so that
$$
\int_{U}\dist(\Gamma(t),\partial\Omega)
^{\tau-1}\,dt=\intl_U g^{\frac{\tau-1}{\alpha-1}}(t)\,dt
=\intl_U g^q(t)\,dt\le\intl_0^L g^q(t)\,dt\le
C_2\|x-y\|^{q\alpha}/L^{q-1}
$$
Since $\|x-y\|\le L$, we obtain
\bel{F-BT} \int_{U}\dist(\Gamma(t),\partial\Omega)
^{\tau-1}\,dt\le C_2\|x-y\|^{q\alpha}/L^{q-1}\le
C_2\|x-y\|^{q\alpha-q+1}= C_2\|x-y\|^{\tau}. \ee
\par Finally, we put
$$\GW:=\Gamma|_U.$$
Then inequalities \rf{PR-C}, \rf{F-B} and \rf{F-BT} shows
that inequalities \rf{H-AA}, \rf{H-LGW} and \rf{H-GW} of
Theorem \reff{S-IM} are satisfied.
\par It remains to prove inequality \rf{H-REG}. We observe
that the set $U=U_k$ is obtained by the standard Cantor
procedure (for $m$-adic intervals; recall that in the
classical case $m=3$). The reader can easily see that this
Cantor set possess the following property: for each
interval $I$ centered in $U$ with $|I|\le 2L$ we have
\bel{R-U}
|I|\le 4m\,|I\cap U|.
\ee
\par Let $B=B(c,r)$ be a ball of radius $r$ centered at a
point $c\in\GW$. We may assume that $r\le\|x-y\|/4$. Then
either $x$ or $y$ does not belong to $B$. Suppose that
$y\notin B$. Then there exists a point $v\in\partial B\cap
\Gamma_{cy}$ such that the arc $\Gamma_{cv}\subset B$.
\par Let us consider the arc $\Gamma_{xc}$\,. If
$\Gamma_{xc}\subset B$, we put
\bel{D1-L} \ell:=\lng(\Gamma_{cv}). \ee
Clearly, in this case $\ell\ge r$.
\par Assume that $\Gamma_{xc}\nsubseteq B$. Then there
exists a point $u\in\partial B\cap \Gamma_{xc}$ such that
the arc $\Gamma_{uc}\subset B$. In this case we put
\bel{D2-L}
\ell:=\min\{\lng(\Gamma_{uc}),\lng(\Gamma_{cv})\}. \ee
Since $\lng(\Gamma_{uc})\ge r,$ $\lng(\Gamma_{cv})\ge r,$
again we have $\ell\ge r$.
\par Recall that $c\in\GW$ so that $c=\Gamma(A)$ for some
$a\in U$. By $I$ we denote the interval
$I:=[a-\ell,a+\ell]$.  Then, by definitions \rf{D1-L},
\rf{D2-L}, the arc $\Gamma(I)\subset B$. Since $I$ is
centered in $U$ and  $|I|=2\ell\le 2L$, by \rf{R-U},
$$ |I|\le 4m\,|I\cap U|. $$
Since $\Gamma(I)\subset B$, we have $|I\cap U|\le
\lng(B\cap\GW)$ so that
$$
\diam B=2r\le 2\ell=|I|\le 4m\,|I\cap U|\le
4m\lng(B\cap\GW)
$$
proving \rf{H-REG}.
\par Theorem \reff{S-IM} is completely proved.\bx
\SECT{3. Sobolev functions on subhyperbolic domains.}{3}
\indent
\par Let us fix some additional notation. In what follows,
the terminology ``cube'' will mean a closed cube in ${\bf
R}^{n}$ whose sides are parallel to the coordinate axes. We
let $Q(x,r)$ denote the cube in $\RN$ centered at $x$ with
side length $2r$. Given $\lambda>0$ and a cube $Q$ we let
$\lambda Q$ denote the dilation of $Q$ with respect to its
center by a factor of $\lambda $. (Thus $\lambda
Q(x,r)=Q(x,\lambda r)$.)
\par It will be convenient for us to measure distances in
$\RN$ in the uniform norm
$$
\|x\|:=\max\{|x_i|:~i=1,...,n\}, \ \ \
x=(x_1,...,x_n)\in\RN.
$$
Thus every cube
$$
Q=Q(x,r):=\{y\in\RN:\|y-x\|\le r\}
$$
is a ``ball" in $\|\cdot\|$-norm  of ``radius" $r$ centered
at $x$.
\par By $\chi_A$ we denote the characteristic function of
$A$; we put $\chi_A\equiv 0$ whenever $A=\emptyset$.
\par Let $\Ac=\{Q\}$ be a family of cubes in $\RN$. By
$M(\Ac)$ we denote its {\it covering multiplicity}, i.e.,
the minimal positive integer $M$ such that every point
$x\in\RN$ is covered by at most $M$ cubes from $\Ac$. Thus
$$
M(\Ac):=\sup_{x\in\RN}\sum_{Q\in \Ac}\chi_Q(x).
$$
\par Recall that $\PK$ denotes the space of polynomials of
degree at most $k$ defined on $\RN$. Also recall that,
given a $k$-times differentiable function $F$ and a point
$x\in\RN,$ we let $T_{x}^{k}(F)$ denote the Taylor
polynomial of $F$ at $x$ of degree at most $k$:
$$ T_{x}^{k}(F)(y):=\sum_{|\beta|\leq
k}\frac{1}{\beta!}(D^{\beta}F)(x)(y-x)^{\beta}~,~~y\in \RN.
$$
\par Let $\Omega$ be a domain in $\RN$ and let $k\in\N$ and
$p\in[1,\infty]$. We let $L_{p}^{k}(\Omega)$ denote the
(homogeneous) Sobolev space of all functions $f\in
L_{1,\,loc}(\Omega)$ whose distributional partial
derivatives on $\Omega$ of order $k$ belong to
$L_{p}(\Omega)$. $L_{p}^{k}(\Omega)$ is normed by
$$
\|f\| _{L_{p}^{k}(\Omega)}:=\left(\intl_{\Omega}\|\nabla^k
f\|^p\,dx\right)^{\frac1p}
$$
where $\nabla^k f$ denotes the vector with components
$D^\beta f$, $|\beta|=k,$ and
$$
\|\nabla^k f\|(x):=\left(\sum_{|\beta|=k}
|D^{\beta }f(x)|^2\right)^{\frac12}\,,
~~~~x\in\Omega.
$$
\par By the Sobolev imbedding theorem, see e.g., \cite{M},
p.\ 60, every $f\in L_{p}^{k}(\Omega),p>n,$ can be
redefined, if necessary, in a set of Lebesgue measure zero
so that it belongs to the space $C^{k-1}(\Omega)$.
Moreover, for every cube $Q\subset\Omega$, every $x,y\in Q$
and every multiindex $\beta,|\beta|\le k-1,$ the following
inequality
\bel{TH-LOC} |D^{\beta}T_{x}^{k-1}(f)(x)-
D^{\beta}T_{y}^{k-1}(f)(x)|\le C(n,p)\|x-y\|
^{k-|\beta|-\tfrac{n}{p}}\left(\,\intl_{Q}\| \nabla^k
f\|^p\,dx\right)^{\frac1p} \ee
holds. In particular, the partial derivatives of order
$k-1$ satisfy a (local) H\"{o}lder condition of order
$\alpha :=1-\frac{n}{p}$:
$$
|D^{\beta }f(x)-D^{\beta }f(y)| \le
C(n,p)\|f\|_{L_{p}^{k}(\Omega)}\|x-y\|^{1-\frac{n}{p}},
~~~|\beta |=k-1,
$$
provided $Q$ is a cube in $\Omega$ and $x,y\in Q$.
\par Thus, for $p>n$, we can identify each element $f\in
L_{p}^{k}(\Omega)$ with its unique $C^{k-1}$-representative
on $\Omega$. This will allow us to restrict our attention
to the case of Sobolev $C^{k-1}$-functions.
\par The main result of this section is the following
\begin{theorem}\lbl{SH-TE} Let $n<p<\infty$,
$\alpha=(p-n)/(p-1),$ and let $\Omega$ be an
$\alpha$-subhyperbolic domain in $\RN$. There exists a
constant $\ps, n<\ps<p,$ and constants $\lambda,\theta,C>0$
depending only on $n,p,k,C_{\alpha,\Omega}$ and
$\theta_{\alpha,\Omega}$, such that the following is true:
Let $f\in L_{p}^{k}(\Omega)$,  $x,y\in\Omega,
\|x-y\|\le\theta,$ and let $Q_{xy}:=Q(x,\|x-y\|)$. Then for
every multiindex $\beta,|\beta|\le k-1,$ the following
inequality
\bel{MR} |D^{\beta}T_{x}^{k-1}(f)(x)-
D^{\beta}T_{y}^{k-1}(f)(x)|\le C\|x-y\|
^{k-|\beta|-\tfrac{n}{\ps}}\left(\, \intl_{(\lambda
Q_{xy})\cap\Omega}\|\nabla^k f\|^\ps\,dx\right)^{\frac1\ps}
\ee
holds.
\end{theorem}
\par {\it Proof.} We will be needed the following
\begin{lemma}\lbl{CH-CR} Let $x,y\in\RN$ and let
$\gamma\subset\Omega$ be a continuous curve joining $x$ to
$y$. There is a finite family of cubes
$\Qc=\{Q_0,...,Q_m\}$ such that:
\par (i). $Q_0\ni x,Q_m\ni y$, $Q_i\ne Q_j$, $i\ne j$,
$0\le i,j\le m$,  and
$$
Q_i\cap Q_{i+1}\ne\emp,~~~~i=0,...,m-1.
$$
\par (ii). For every cube $Q=Q(z,r)\in\Qc$ we have
$z\in\gamma$ and $r=\tfrac18\dist(z,\partial\Omega)$.
\par (iii). For each $Q\in\Qc$ the cube $2Q\subset\Omega$.
Moreover, the covering multiplicity of the family of cubes
$ 2\Qc:=\{2Q:~Q\in\Qc\} $ is bounded by a constant
$C=C(n)$.
\end{lemma}
\par {\it Proof.} For every $z\in\Gamma$ we let $Q^{(z)}$
denote the cube
$$
Q^{(z)}:=Q(z,\tfrac18\dist(z,\partial\Omega)).
$$
We put $\Ac:=\{Q^{(z)}:~z\in\Gamma\}$. By the Besicovitch
covering theorem, see e.g. \cite{G}, there exists a finite
subcollection $\Bc\subset\Ac$ such that $\Bc$ still covers
$\Gamma$ but no point which lies in more than $C(n)$ of the
cubes of $\Bc$. (Thus the covering multiplicity $M(\Bc)\le
C(n).$)
\par Given $Q',Q''\in\Bc$ we write $Q'\sim Q''$ if there
exists a family $\{K_0,...,K_\ell\}\subset\Bc$ such that
$K_0=Q'$, $K_\ell=Q''$, $K_i\ne K_j, i\ne j,$ $0\le i,j\le
\ell$, and $K_i\cap K_{i+1}\ne\emp$ for every
$i=0,...,\ell-1$. Fix a cube $\tQ\in\Bc$ such that
$x\in\tQ,$ and put
\bel{D-AC} \Bc':=\{Q\in\Bc:~Q\sim\tQ\} \ee
and $E:=\cup\{Q:~Q\in\Bc'\}.$
\par Prove that $y\in E$. Assume that this is not true,
i.e., $y\notin E.$ Since $\tQ\ni x$, we have $x\in E$ so
that $E\cap \gamma\ne\emp$. Put
$F:=\cup\{Q:~Q\in\Bc\setminus \Bc'\}.$ Observe that every
cube $Q\in\Bc$ such that $Q\cap E\ne\emp$ belongs to $\Bc'$
so that $E\cap F=\emp$.
\par On the other hand, since $\gamma\subset E\cup F$ and
$y\notin E$, we have $y\in F$ proving that
$F\cap\gamma\ne\emp.$ Thus the sets $E\cap\gamma$ and
$F\cap\gamma$ is a partition of the continuous curve
$\gamma$ into two closed disjoint sets; a contradiction.
We have proved that $y\in Q^*$ for some $Q^*\in\Bc'$ so
that, by definition \rf{D-AC}, there exists a family of
cubes $\Qc=\{Q_0,,...,Q_m\}$ satisfying conditions (i) and
(ii) of the the lemma.
\par Prove (iii). Let $Q=Q(z,r)\in\Qc$. Then $z\in\gamma$
and $r=\tfrac18\dist(z,\partial\Omega)$ so that for every
$u\in 2Q=Q(z,2r)$ we have
$$
\dist(z,\partial\Omega)\le
\dist(u,\partial\Omega)+\|u-z\|\le
\dist(u,\partial\Omega)+2r\le
\dist(u,\partial\Omega)+\tfrac14\dist(z,\partial\Omega).
$$
Hence,
$$
0<\tfrac34\dist(z,\partial\Omega)\le
\dist(u,\partial\Omega)
$$
proving that $2Q\subset\Omega$.
\par It remains to prove that the covering multiplicity
$M(2\Qc)\le C(n)$. We know that $M(\Qc)\le M(\Bc)\le C(n)$.
Fix a cube $Q=Q(z,r)\in\Qc$. Let $Q_i=Q(z_i,r_i)\in\Qc$ be
an arbitrary cube such that
\bel{I-Q1} (2Q)\cap (2Q_i)\ne\emp. \ee
Then $\|z-z_i\|\le 2r+2r_i$ so that
$$
r=\tfrac18\dist(z,\partial\Omega)\le
\tfrac18\dist(z_i,\partial\Omega)+\tfrac18\|z-z_i\|\le
r_i+\tfrac18(2r+2r_i)=\tfrac14 r+\tfrac54 r_i.
$$
Hence $r\le\tfrac53r_i$. In the same way we prove that
$r_i\le\tfrac53r$.
\par Since $\Qc$ has the covering multiplicity at most
$C(n)$, this collection of cubes can be partitioned into at
most $N(n)$ families of pairwise disjoint cubes, see e.g.
\cite{BrK}. Therefore, without loss of generality, we may
assume that $\Qc$ itself is a collection of pairwise
disjoint cubes.
\par Since $\tfrac12r\le r_i\le 2r$, we have $2^{-n}|Q|\le
|Q_i|\le 2^{n}|Q|$.  Also, by \rf{I-Q1} and the inequality
$r_i\le 2r$, we have $Q_i\subset 7Q$. Thus the number of
cubes $Q_i$ satisfying \rf{I-Q1} is bounded by
$|7Q|/(2^{-n}|Q|)=2^n7^n$. The lemma is proved.\bx
\par Let $x,y\in\Omega,$ $\|x-y\|\le \theta,$ where
$\theta$ is the constant from Theorem \reff{S-IM}. By this
theorem there exist constants $\as=\as(n,p),0<\as<\alpha,$
and $C=C(n,p)>0$ such that for every $\ve>0$ there exists a
rectifiable curve $\Gamma\subset\Omega$ and a finite family
of arcs $\GW\subset\Gamma$ satisfying conditions (i),(ii)
of the theorem.
\par Observe that, by inequality \rf{H-GW} (with
$\tau=\alpha$) and by \rf{H-AA},
\bel{E-AL} \intl_{\Gamma} \dist(z,\partial\Omega)
^{\alpha-1}\,ds(z)\le C\|x-y\|^{\alpha}. \ee
Also, by part (ii) of Theorem \reff{S-IM},
\bel{LG-1}\lng(\Gamma)\le \tC\|x-y\|\ee
with $\tC=2e^{2C_{\alpha,\Omega}},$ see \rf{LGM}.
\par By Lemma \reff{CH-CR}, there exists a collection of
cubes
$$
\Qc=\{Q_0,...,Q_m\}
$$
satisfying conditions  (i)-(iii) of the lemma. In
particular, by (i), $Q_0\ni x, Q_m\ni y,$ $Q_i\ne Q_j, i\ne
j,$ and
$$
Q_i\cap Q_{i+1}\ne\emp,~~~i=0,...,m-1.
$$
Let $Q_i=Q(z_i,r_i), i=0,...,m,$ (recall that by (ii)
$z_i\in\Gamma$ and
$r_i=\tfrac18\dist(z_i,\partial\Omega)$). Let $a_i\in
Q_{i-1}\cap Q_{i},$ $i=1,...,m.$ We also put
$a_0:=x,a_{m+1}:=y$.
\par We may assume that for every $Q\in\Qc$ either $x\notin
Q$ or $y\notin Q$. In fact, otherwise $x,y\in Q$. But by
condition (iii) of Lemma \reff{CH-CR}, $2Q\subset \Omega$.
Then the cube $Q(x,\|x-y\|)\subset 2Q\subset\Omega$ as
well. It remains to apply inequality \rf{TH-LOC} to $x$ and
$y$ with $p=\ps, n<\ps<p,$ and \rf{MR} follows.
\par Thus we assume that for every cube
$Q_i=Q(z_i,r_i)\in\Qc$
$$
\text{either }~~x\notin Q_i~~ \text{~or }~~y\notin Q_i.
$$
Since $x,y,z_i\in \Gamma,$ we have $\partial Q_i\cap
\Gamma\ne\emp$ so that there exists a point $a_i\in\partial
Q_i\cap \Gamma$. Hence, by \rf{LG-1},
\bel{LA} r_i=\|z_i-a_i\|\le\lng(\Gamma\cap Q_i) \ee
so that, by \rf{LG-1},
\bel{R-XY} r_i\le\lng (\Gamma)\le \tC\|x-y\|,~~~i=0,...,m.
\ee
\par Now we have
\be A&:=&|D^{\beta}T_{x}^{k-1}(f)(x)-
D^{\beta}T_{y}^{k-1}(f)(x)|=
|D^{\beta}T_{a_0}^{k-1}(f)(a_0)-
D^{\beta}T_{a_{m+1}}^{k-1}(f)(a_0)|\nn\\& \le& \sum_{i=0}^m
|D^{\beta}T_{a_i}^{k-1}(f)(a_0)-
D^{\beta}T_{a_{i+1}}^{k-1}(f)(a_0)|. \nn\ee
Put
$$
P_i(z):=T_{a_i}^{k-1}(f)(z)-
T_{a_{i+1}}^{k-1}(f)(z),~~~i=0,...,m.
$$
The polynomial $P_i\in\Pc_{k-1}$ so that for every
multiindex $\beta,|\beta|\le k-1,$ we have
$$
D^\beta P_i(z)=\sum_{|\eta|\le\, k-1-|\beta|}
\frac{1}{\eta!}\,
D^{\eta+\beta} P_i(a_i)\,(z-a_i)^\eta,~~~z\in\RN.
$$
Hence
$$
|D^\beta P_i(a_0)|\le C\sum_{|\eta|\le\, k-1-|\beta|}
|D^{\eta+\beta} P_i(a_i)|\,\|a_0-a_i\|^{|\eta|}.
$$
We put
\bel{Q-1} \Qc_1:=\{Q\in\Qc:~Q\cap\GW\ne\emp\},~~~~~
I_1:=\{i\in\{0,...,m\}:~Q_i\in\Qc_1\}, \ee
and
\bel{Q-2} \Qc_2:=\Qc\setminus \Qc_1,~~~~~
I_2:=\{i\in\{0,...,m\}:~Q_i\in\Qc_2\}. \ee
Then
\be
A&\le& \sum_{i=0}^m  |D^\beta P_i(a_0)|\le
C\sum_{i=0}^m \sum_{|\eta|\le\, k-1-|\beta|}
|D^{\eta+\beta} P_i(a_i)|\,\|a_0-a_i\|^{|\eta|}\nn\\
&=& C\sum_{|\eta|\le\, k-1-|\beta|}\left(\sum_{i=0}^m
|D^{\eta+\beta} P_i(a_i)|\,\|a_0-a_i\|^{|\eta|}\right).
\nn
\ee
Let
$$
A'_\eta:=\sum_{i\in I_1} |D^{\eta+\beta}
P_i(a_i)|\,\|a_0-a_i\|^{|\eta|}
$$
and
$$
A''_\eta:=\sum_{i\in I_2} |D^{\eta+\beta}
P_i(a_i)|\,\|a_0-a_i\|^{|\eta|}.
$$
We have proved that
\bel{E-AM}
A\le C\sum_{|\eta|\le\, k-1-|\beta|}
(A'_\eta+A''_\eta).
\ee
\par Let $|\eta|\le\, k-1-|\beta|$. Our next aim is to show
that
\bel{E-A1}
A'_\eta\le C\,\|x-y\|^{k-|\beta|-\frac{n}{\ps}}
\left(\,\, \intl_{(\lambda Q_{xy})\cap\Omega} \|\nabla^k
f\|^\ps\,dx\right)^{\frac1\ps}
\ee
where
$$\ps:=\frac{n-\as}{1-\as}$$
and $\lambda:=2\tC$. (Recall that $Q_{xy}:=Q(x,\|x-y\|).$)
Also we will prove that
\bel{E-A2} A''_\eta\le C \|x-y\|^{k-|\beta|-\frac{n}{p}}
\left(\,\,\intl_{E} \|\nabla^kf\|^p\,dx\right)^{\frac1p}
\ee
where is $E$ is a subset of $\Omega$ of the Lebesgue
measure  $|E|\le C\ve^n$.
\par Since
$$
0<\as=\frac{\ps-n}{\ps-1}<\alpha=\frac{p-n}{p-1}<1,
$$
we have $n<\ps<p$. Since $a_i,a_{i+1}\in Q_i=Q(z_i,r_i)$,
$i\in I_1$, by inequality \rf{TH-LOC} (with $p=\ps$),
\be |D^{\eta+\beta}P_i(a_i)|&=&|D^{\eta+\beta}
T_{a_i}^{k-1}(f)(a_i)-
D^{\eta+\beta}T_{a_{i+1}}^{k-1}(f)(a_i)| \nn\\&\le&
C\,\|a_i-a_{i+1}\|^{k-|\eta|-|\beta|-\frac{n}{\ps}}
\left(\,\,\intl_{Q_i}\|\nabla^k
f\|^\ps\,dx\right)^{\frac1\ps}\nn\ee
so that
\bel{AP} |D^{\eta+\beta}P_i(a_i)|\le
C\,r_i^{k-|\eta|-|\beta|-\frac{n}{\ps}}
\left(\,\,\intl_{Q_i}\|\nabla^k
f\|^\ps\,dx\right)^{\frac1\ps},~~~i\in I_1. \ee
In a similar way we prove that
$$
|D^{\eta+\beta}P_i(a_i)|\le
C\,r_i^{k-|\eta|-|\beta|-\frac{n}{p}}
\left(\,\,\intl_{Q_i}\|\nabla^k
f\|^p\,dx\right)^{\frac1p},~~~i\in I_2.
$$
\par Prove that for each $i\in\{0,...,m\}$ we have
\bel{A-A}
\|a_0-a_i\|\le 2\tC\|x-y\|.
\ee
In fact,
$$
\|a_0-a_i\|=\|x-a_i\|\le
\|x-z_i\|+\|z_i-a_i\|\le\|x-z_i\|+r_i.
$$
Since $x,z_i\in\Gamma,$ by \rf{LG-1},
$\|x-z_i\|\le\lng(\Gamma)\le\tC\|x-y\|$. Also, by
\rf{R-XY}, $r_i\le\tC\|x-y\|$, proving \rf{A-A}.
\par Hence, by \rf{AP},
$$
A'_\eta\le (2\tC)^{|\eta|}\sum_{i\in I_1} |D^{\eta+\beta}
P_i(a_i)|\|x-y\|^{|\eta|}\le C\sum_{i\in I_1}
r_i^{k-|\eta|-|\beta|-\frac{n}{\ps}}\|x-y\|^{|\eta|}
\left(\,\,\intl_{Q_i}\|\nabla^k
f\|^\ps\,dx\right)^{\frac1\ps}.
$$
Since $k-|\eta|-|\beta|\ge 1$ and $r_i\le\tC\|x-y\|$, we
have
\be r_i^{k-|\eta|-|\beta|-\frac{n}{\ps}}\|x-y\|^{|\eta|}
&=&r_i^{k-|\eta|-|\beta|-1}r_i^{1-\frac{n}{\ps}}
\|x-y\|^{|\eta|}\nn\\ &\le& C\,r_i^{1-\frac{n}{\ps}}
\|x-y\|^{k-|\eta|-|\beta|-1}
\|x-y\|^{|\eta|}\nn\\&=&C\,r_i^{1-\frac{n}{\ps}}
\|x-y\|^{k-|\beta|-1}.\nn\ee
Hence
$$
A'_\eta\le C\,\|x-y\|^{k-|\beta|-1}
\sum_{i\in I_1} r_i^{1-\frac{n}{\ps}}
\left(\,\,\intl_{Q_i}\|\nabla^k
f\|^\ps\,dx\right)^{\frac1\ps}.
$$
By the H\"{o}lder inequality,
$$
A'_\eta\le C\,\|x-y\|^{k-|\beta|-1}
 \left(\,\sum_{i\in I_1}(r_i^{1-\frac{n}{\ps}})^{\frac{\ps}{\ps-1}}
\right)^{1-\frac1\ps}\left(\sum_{i\in I_1}
\intl_{Q_i}\|\nabla^k
f\|^\ps\,dx\right)^{\frac1\ps}
$$
so that
$$ A'_\eta\le C\,\|x-y\|^{k-|\beta|-1} \left(\,\sum_{i\in
I_1} r_i^{\as}
\right)^{1-\frac1\ps}\left(M(\Qc_1)\intl_{U_1}\|\nabla^k
f\|^\ps\,dx\right)^{\frac1\ps} $$
where
$$
U_1:=\cup \{Q:~Q\in\Qc_1\}.
$$
Recall that $M(\Qc_1)$ stands for the covering multiplicity
of the collection $\Qc_1$. Since $M(\Qc_1)\le M(\Qc)\le
C(n),$ we obtain
\bel{W-A1} A'_\eta\le C\,\|x-y\|^{k-|\beta|-1}
\left(\,\sum_{i\in I_1} r_i^{\as}
\right)^{1-\frac1\ps}\left(\,\,\intl_{U_1}\|\nabla^k
f\|^\ps\,dx\right)^{\frac1\ps} \ee
In a similar way we prove that
\bel{W-A2} A''_\eta\le C\,\|x-y\|^{k-|\beta|-1}
\left(\,\sum_{i\in I_2} r_i^{\alpha}
\right)^{1-\frac1p}\left(\,\,\intl_{U_2}\|\nabla^k
f\|^p\,dx\right)^{\frac1p} \ee
where
$$
U_2:=\cup \{Q:~Q\in\Qc_2\}.
$$
\par Let us prove that
\bel{R-AS}
\sum_{i\in I_1} r_i^{\as}\le C\,\|x-y\|^{\as}
\ee
and
\bel{R-AS2}
\sum_{i\in I_2} r_i^{\alpha}\le
C\,\|x-y\|^{\alpha}.
\ee
\par We begin with the proof of inequality \rf{R-AS}. Let
$i\in I_1$ and let
$$
Q_i=Q(z_i,r_i),~~~r_i=\tfrac18\dist(z_i,\partial\Omega).
$$
Recall that $Q_i\cap\GW\ne\emp$ so that there exist a point
$b_i\in Q_i\cap\GW$. By \rf{H-REG},
$$
r_i\le C\,\lng(\GW\cap Q(b_i,r_i)).
$$
But $Q(b_i,r_i)\subset 2Q_i=Q(z_i,2r_i)$ so that
$$
r_i\le C\,\lng(2Q_i\cap\GW ).
$$
Hence,
$$
r_i^{\as}\le C\,r_i^{\as-1}\lng(2Q_i\cap\GW).
$$
\par Let $z\in 2Q_i\cap\GW$. Then
$$
|\dist(z,\partial\Omega)-\dist(z_i,\partial\Omega)|
\le\|z-z_i\|\le 2r_i.
$$
Since $\dist(z_i,\partial\Omega)=8r_i$, we have
$$
\dist(z,\partial\Omega)\le \dist(z_i,\partial\Omega)
+2r_i= 10r_i.
$$
Hence,
$$
r_i^{\as-1}\le C\dist(z,\partial\Omega)^{\as-1},
~~~~z\in 2Q_i\cap\GW,
$$
so that
$$
r_i^{\as}\le C\intl_{2Q_i\cap\GW}
\dist(z,\partial\Omega)^{\as-1}\,ds(z).
$$
We put $2\Qc_1:=\{2Q:~Q\in\Qc_1\}$. We have
$$
\sum_{i\in I_1} r_i^{\as}\le C\sum_{i\in
I_1}\intl_{2Q_i\cap\GW}
\dist(z,\partial\Omega)^{\as-1}\,ds(z)
$$
so that
$$
\sum_{i\in I_1} r_i^{\as}\le C\,M(2\Qc_1) \intl_{\GW}
\dist(z,\partial\Omega)^{\as-1}\,ds(z).
$$
\par By part (iii) of Lemma \reff{CH-CR},  the covering
multiplicity
$$
M(2\Qc_1)\le M(2\Qc)\le C(n).
$$
Hence,
$$
\sum_{i\in I_1} r_i^{\as}\le C \intl_{\GW}
\dist(z,\partial\Omega)^{\as-1}\,ds(z).
$$
\par In a similar way we prove that
$$
\sum_{i\in I_2} r_i^{\alpha}\le C \intl_{\Gamma}
\dist(z,\partial\Omega)^{\alpha-1}\,ds(z).
$$
Combining these inequalities with \rf{H-GW} (where we put
$\tau=\as$) and \rf{E-AL}, we obtain the required
inequalities \rf{R-AS} and \rf{R-AS2}.
\par Hence, by \rf{W-A1},
\be A'_\eta&\le& C\,\|x-y\|^{k-|\beta|-1}
(\|x-y\|^{\as})^{1-\frac1\ps}\left(\,\,\intl_{U_1}\|\nabla^k
f\|^\ps\,dx\right)^{\frac1\ps}\nn\\
&=& C\,\|x-y\|^{k-|\beta|-\frac{n}{\ps}}
\left(\,\,\intl_{U_1}\|\nabla^k
f\|^\ps\,dx\right)^{\frac1\ps} \nn\ee
and, by \rf{W-A2},
\bel{V-2} A''_\eta\le C\,\|x-y\|^{k-|\beta|-\frac{n}{p}}
\left(\,\,\intl_{U_2}\|\nabla^k f\|^p\,dx\right)^{\frac1p}.
\ee
\par Recall that for each $Q=Q(z,r)\in\Qc$ its center, the
point $z$, belongs to $\Gamma$. Moreover, by \rf{R-XY},
$r\le\tC\|x-y\|$, and, by \rf{LG-1}, $\lng
(\Gamma)\le\tC\|x-y\|$. Hence,
$$
Q\subset Q(x,2\tC\|x-y\|)=\lambda Q_{xy},~~~~Q\in\Qc,
$$
with $\lambda:=2\tC$. Now we have
$$
U_1:=\cup \{Q:~Q\in\Qc_1\}\subset (\lambda Q_{xy})\cap\Omega
$$
so that
$$
A'_\eta\le C\,\|x-y\|^{k-|\beta|-\frac{n}{\ps}}
\left(\,\,
\intl_{(\lambda Q_{xy})\cap\Omega}
\|\nabla^k f\|^\ps\,dx\right)^{\frac1\ps}
$$
proving \rf{E-A1}.
\par Let us put $E:=U_2$ and prove that $|E|\le C\,\ve^n$.
We have
\bel{M-E}
|E|=|\cup\{Q:~Q\in\Qc_2\}|\le
\sum\{|Q|:~Q\in\Qc_2\} =2^n\sum_{I\in I_2}r_i^n \le
2^n\left(\sum_{I\in I_2}r_i\right)^n.
\ee
\par By \rf{Q-1} and \rf{Q-2},
$$
Q_i\cap\GW=\emp~~\text{for every}~~ Q_i\in\Qc_2,
$$
so that
$$
Q_i\cap\Gamma=Q_i\cap(\Gamma\setminus\GW).
$$
But, by \rf{LA}, $ r_i\le\lng(Q_i\cap\Gamma) $ so that
$r_i\le\lng(Q_i\cap(\Gamma\setminus\GW))$. Hence
$$
\sum_{I\in I_2}r_i\le \sum_{I\in I_2}
\lng(Q_i\cap(\Gamma\setminus\GW))\le M(\Qc_2)
\lng(\Gamma\setminus\GW).
$$
Since $M(\Qc_2)\le M(\Qc)\le C(n),$ we obtain
$$
\sum_{I\in I_2}r_i\le  C\,
\lng(\Gamma\setminus\GW)
$$
so that, by \rf{H-LGW},
$$
\sum_{I\in I_2}r_i\le C\,\ve.
$$
Combining this inequality with \rf{M-E}, we obtain the
required inequality $|E|\le C\,\ve^n$. This inequality and
\rf{V-2} imply \rf{E-A2}.
\par Now, by \rf{E-AM}, \rf{E-A1} and \rf{E-A2},
\be
A&:=&|D^{\beta}T_{x}^{k-1}(f)(x)-
D^{\beta}T_{y}^{k-1}(f)(x)|\le C\sum_{|\eta|\le\,
k-1-|\beta|} (A'_\eta+A''_\eta)\nn\\
&\le& C\sum_{|\eta|\le\, k-1-|\beta|}
\|x-y\|^{k-|\beta|-\frac{n}{\ps}} \left(\,\,
\intl_{(\lambda Q_{xy})\cap\Omega}
\|\nabla^k f\|^\ps\,dx\right)^{\frac1\ps}\nn\\
&+& C\sum_{|\eta|\le\, k-1-|\beta|}
\|x-y\|^{k-|\beta|-\frac{n}{p}}
\left(\,\,\intl_{E}\|\nabla^kf\|^p\,dx\right)^{\frac1p}.
\nn
\ee
We obtain
\be
A&\le& C \|x-y\|^{k-|\beta|-\frac{n}{\ps}} \left(\,\,
\intl_{(\lambda Q_{xy})\cap\Omega}
\|\nabla^k f\|^\ps\,dx\right)^{\frac1\ps}\nn\\
&+& C\|x-y\|^{k-|\beta|-\frac{n}{p}}
\left(\,\,\intl_{E}\|\nabla^kf\|^p\,dx\right)^{\frac1p}.
\nn \ee
But $f\in L_{p}^{k}(\Omega)$ so that
$$
\intl_{\Omega}\|\nabla^kf\|^p\,dx=\|f\|^p
_{L_{p}^{k}(\Omega)}<\infty.
$$
Hence,
$$
\intl_{E}\|\nabla^kf\|^p\,dx\to 0~~\text{as}~~~
|E|=C\ve^n\to 0,
$$
proving the theorem.\bx
\SECT{4. Extension of Sobolev functions defined on
subhyperbolic domains.}{4}
\indent
\par Given a cube $Q\subset\RN$ and a function $f\in
L_{q}(Q)$, $0<q\le\infty,$ we let $\Ec_k(f;Q)_{L_q}$ denote
the {\it normalized local best approximation} of $f$ on $Q$
in $L_q$-norm by polynomials of degree at most $k-1$, see
Brudnyi \cite{Br1}. More explicitly, we define
$$
\Ec_k(f;Q)_{L_q}:=|Q|^{-\frac{1}{q}}\inf_{P\in\Pc_{k-1}}
\|f-P\|_{L_q(Q)}=\inf_{P\in\Pc_{k-1}}\left(\frac{1}{|Q|}
\int_Q|f-P|^q\,dx\right)^{\frac{1}{q}}.
$$
\par In the literature $\Ec_k(f;Q)_{L_q}$ is also sometimes
called the {\it local oscillation} of $f$, see e.g. Triebel
\cite{T2}.
\par Given a locally integrable function $f$ on $\RN$, we
define its {\it sharp maximal function} $f^{\sharp}_k$ by
letting
$$
f^{\sharp}_k(x):=\sup_{r>0}
r^{-k}\,\Ec_{k}(f;Q(x,r))_{L_1}.
$$
\par Recall that a function $f\in\WKP$, $1<p\le \infty$, if
and only if $f$ and $f^{\sharp}_k$ are both in $L_p(\RN)$,
see Calder\'{o}n \cite{C}. Moreover, up to constants
depending only on $n,k$ and $p$ the following equivalence,
\bel{AN} \|f\|_{\WKP}\sim
\|f\|_{L_p(\RN)}+\|f_{k}^{\sharp}\|_{L_p(\RN)}, \ee
holds.
\par This characterization motivates the following
definition. Let $S$ be a measurable subset of $\RN$. Given
a function $f\in L_{q,\,loc}(S),$ and a cube $Q$ whose
center is in $S$, we let $\Ec_k(f;Q)_{L_q(S)}$ denote the
norma\-lized best approximation of $f$ on $Q$ in
$L_q(S)$-norm:
\bel{ES} \Ec_k(f;Q)_{L_q(S)}:=|Q|^{-\frac{1}{q}}
\inf_{P\in\PK}\|f-P\|_{L_q(Q\cap S)}
=\inf_{P\in\Pc_{k-1}}\left(\frac{1}{|Q|} \int_{Q\cap S}
|f-P|^q\,dx\right)^{\frac{1}{q}}. \ee
By $f^{\sharp}_{k,S},$ we denote the sharp maximal function
of $f$ on $S$,
$$ f^{\sharp}_{k,S}(x):=\sup_{r>0}
r^{-k}\,\Ec_{k}(f;Q(x,r))_{L_1(S)},\ \ \ \ \ \ x\in S.
$$
(Thus, $f^{\sharp}_{k}=f^{\sharp}_{\alpha,\RN}$~.)
\par Let $\Omega\subset\RN$ be a subhyperbolic domain. The
following two corollaries of Theorem \reff{SH-TE} present
estimates of the local best approximations and the sharp
maximal function of a function $f\in\WKPO$ via the local
$L_p$-norms and the maximal function of $\nabla^kf$.
\begin{corollary}\lbl{COR1} Let $n<p<\infty$,
$\alpha=(p-n)/(p-1),$ and let $\Omega$ be an
$\alpha$-subhyperbolic domain in $\RN$. There exists a
constant $p^*\in(n,p)$ and constants $\theta,\lambda,C>0$
depending only on $n,p,k,C_{\alpha,\Omega}$ and
$\theta_{\alpha,\Omega}$, such that the following is true:
Let $f\in L_{p}^{k}(\Omega)$. Then for every cube
$Q=Q(x,r)$ with  $x\in\Omega$ and $0<r\le\theta$ the
following inequality
$$
r^{-k}\Ec_k(f;Q)_{L_\infty(\Omega)}
\le C\left(\frac{1}{|\lambda Q|}
\intl_{(\lambda Q)\cap\Omega}\|\nabla^k f\|^{p^*}\,
dx\right)^{\frac{1}{p^*}}
$$
holds.
\end{corollary}
\par {\it Proof.} Let $p^*$ and $\theta$ be the constant
from Theorem \reff{SH-TE}. Let $y\in Q(x,r)$ so that
$\|y-x\|\le r\le \theta$. Applying Theorem \reff{SH-TE} to
the points $y,x$ (with $\beta=0$), we obtain
$$
|f(y)- T_{x}^{k-1}(f)(y)|\le
C\|x-y\|^{k-\frac{n}{p^*}} \left(\,\,\intl_{(\lambda
Q_{xy})\cap\Omega}\|\nabla^k f\|^{p^*}\,dx\right)
^{\frac{1}{p^*}}.
$$
Recall that $Q_{xy}:=Q(x,\|x-y\|)$.
\par Since $n<p^*$, we have $k-\frac{n}{p^*}>0$ so that
\be
|f(y)- T_{x}^{k-1}(f)(y)|&\le& Cr^{k-\frac{n}{p^*}}
\left(\,\, \intl_{(\lambda Q_{xy})\cap\Omega} \|\nabla^k
f\|^{p^*}\,dx\right)^{\frac{1}{p^*}}\nn\\ &\le& Cr^{k}
\left(\,\,\frac{1}{|\lambda Q|} \intl_{(\lambda Q)
\cap\Omega} \|\nabla^k
f\|^{p^*}\,dx\right)^{\frac{1}{p^*}}.
\nn\ee
Hence,
\be
\Ec_k(f;Q)_{L_\infty(\Omega)}&:=&\inf_{P\in\Pc_{k-1}}
\sup_{y\in Q\cap\Omega}|f(y)-P(y)|\le \sup_{y\in
Q\cap\Omega}|f(y)- T_{x}^{k-1}(f)(y)|\nn\\
&\le&  Cr^{k} \left(\,\,\frac{1}{|\lambda Q|}
\intl_{(\lambda Q)\cap\Omega} \|\nabla^k
f\|^{p^*}\,dx\right)^{\frac{1}{p^*}} \nn \ee
proving the corollary.\bx
\par Given a function $g$ defined on $\Omega$ we let
$g^\cw$ denote its extension by zero to all of $\RN$. Thus
$g^\cw(x):=g(x), x\in\Omega,$ and $g^\cw(x):=0,
x\notin\Omega.$
\par As usual, given a function $u\in L_{1,loc}(\RN)$ by
$\Mc f$ we denote the Hardy-Littlewood maximal function
$$
\Mc[f](x):=\sup_{t>0}\frac{1}{|Q(x,t)|}
\int_{Q(x,t)}|f(y)|dy.
$$
\begin{corollary}\lbl{COR2} Let $n<p<\infty$,
$\alpha=(p-n)/(p-1),$ and let $\Omega$ be an
$\alpha$-subhyperbolic domain in $\RN$. There exists a
constant $p^*, n<p^*<p,$ such that for every function $f\in
L_{p}^{k}(\Omega)$ and every $x\in\Omega$ the following
inequality
$$
f^{\sharp}_{k,\Omega}(x)\le
C\left\{(\Mc[(\|\nabla^kf\|^\cw)^{p^*}](x))^{\frac{1}{p^*}}+
\Mc[f^\cw](x)\right\}
$$
holds. The constants $p^*$ and $C$ depend only on
$n,p,k,C_{\alpha,\Omega}$ and $\theta_{\alpha,\Omega}$.
\end{corollary}
\par {\it Proof.} Let $p^*,\lambda$ and $\theta$ be the
constants from Corollary \reff{COR1}. By this corollary,
\be
\sup_{0<r\le\theta}r^{-k}\Ec_k(f;Q(x,r))_{L_1(\Omega)}
&\le&\sup_{0<r\le\theta}
r^{-k}\Ec_k(f;Q(x,r))_{L_\infty(\Omega)}
\nn\\
&\le& C\sup_{0<r\le\theta}\left(\frac{1}{|Q(x,\lambda r)|}
\intl_{Q(x,\lambda r)\cap\Omega}\| \nabla^k f\|^{p^*}
\,dx\right)^{\frac1{p^*}}\nn\\
&\le& C \left\{(\Mc[(\|\nabla^kf\|^\cw)^{p^*}](x))
^{\frac{1}{p^*}}\right\}. \nn\ee
\par On the other hand, by \rf{ES},
\be
\sup_{r>\theta}r^{-k}\Ec_k(f;Q(x,r))_{L_1(\Omega)}
&\le&\theta^{-k}
\sup_{r>\theta}\Ec_k(f;Q(x,r))_{L_1(\Omega)}
\nn\\
&\le& \theta^{-k}\sup_{r>\theta}\left(\frac{1}{|Q(x,r)|}
\intl_{Q(x,r)\cap\Omega}|f|\,dx\right)\nn\\
&\le&  \theta^{-k}\Mc(f^\cw)(x)
\nn\ee
proving the lemma.\bx
\par In \cite{S4} we show that the restrictions of Sobolev
functions to {\it regular subsets} of $\RN$ can be
described in a way similar to the Calder\'{o}n's criterion
\rf{AN}, i.e., via $L_p$-norms of a function and its sharp
maximal function on a set. We recall that a measurable set
$S\subset\RN$ is said to be {\it regular} if there are
constants $\sigma_S\ge 1$ and $\delta_S>0$ such that, for
every cube $Q$ with center in $S$ and with diameter $\diam
Q\le\delta_S$,
$$
|Q|\le \sigma_S |Q\cap S|.
$$
\begin{theorem}\label{EXT1}(\cite{S4}) Let $S$ be a regular
subset of $\RN$. Then a function $f\in L_p(S),~1<p\le
\infty,$ can be extended to a function $F\in \WKP$ if and
only if its sharp maximal function $f_{k,S}^{\sharp}\in
L_p(S).$
In addition,
$$
\|f\|_{\WKP|_S}\sim
\|f\|_{L_p(S)}+\|f_{k,S}^{\sharp}\|_{L_p(S)}
$$
with constants of equivalence depending only on
$n,k,p,\sigma_S$ and $\delta_S$.
\end{theorem}
\medskip
\par {\it Proof of Theorem \reff{MAIN-EXT}.} By inequality
\rf{F-IN}, $\Omega$ is an $\alpha$-subhyperbolic domain
with $\alpha=\frac{p-n}{p-1},$ so that, by Lemma \reff{BP},
$\Omega$ is a {\it regular} subset of $\RN$.
\par Let $p^*\in(n,p)$ be the constant from Corollary
\reff{COR2}. Let $q>p^*$ and let $f\in\WKQO$. We have to
prove that $f$ can be extended to a function $F\in\WKQ$.
Since $\Omega$ is regular, by Theorem \reff{EXT1} it
suffices to show that the sharp maximal function $
f^{\sharp}_{k,\Omega}$ belongs to $L_q(\Omega)$.
\par By Corollary \reff{COR2},
$$
f^{\sharp}_{k,\Omega}(x)\le
C\left\{(\Mc[(\|\nabla^kf\|^\cw)^{p^*}](x))^{\frac{1}{p^*}}+
\Mc[f^\cw](x)\right\},~~~~ x\in\Omega,
$$
so that
\be \|f^{\sharp}_{k,\Omega}\|_{L_q(\Omega)}&\le&
C\left\{\|(\Mc[(\|\nabla^kf\|^\cw)^{p^*}])
^{\frac{1}{p^*}}\|_{L_q(\Omega)}+
\|\Mc[f^\cw]\|_{L_q(\Omega)}\right\}\nn\\
&\le& C\left\{\|(\Mc[(\|\nabla^kf\|^\cw)^{p^*}])
^{\frac{1}{p^*}}\|_{L_q(\RN)}+
\|\Mc[f^\cw]\|_{L_q(\RN)}\right\}. \nn \ee
By the Hardy-Littlewood maximal theorem
$$
\|\Mc[f^\cw]\|_{L_q(\RN)}\le
C\|f^\cw\|_{L_q(\RN)}=C\|f\|_{L_q(\Omega)}.
$$
(Recall that $f^\cw$ denotes the extension of $f$ by zero
to all of $\RN$.)
\par Applying the Hardy-Littlewood maximal theorem to the
function $g:=(\|\nabla^kf\|^\cw)^{p^*}$ in the space
$L_s(\RN)$ with $s:=q/p^*>1$, we obtain
\be A:=\|(\Mc[(\|\nabla^kf\|^\cw)^{p^*}])
^{\frac{1}{p^*}}\|_{L_q(\RN)}&=&
\left(\|\Mc[g]\|_{L_s(\RN)}\right)^{\frac{1}{p^*}} \le
C\left(\|g\|_{L_s(\RN)}\right)^{\frac{1}{p^*}}\nn\\
&=& C\left\{\left(\,\,\intl_{\RN}
\left[\left(\|\nabla^kf\|^\cw\right)^{p^*}\right]
^{\frac{q}{p^*}}\,dx\right)^{\frac{p^*}{q}}
\right\}^{\frac{1}{p^*}} \nn\ee
so that
$$
A\le C\left(\,\,\intl_{\RN}
\left(\|\nabla^kf\|^\cw\right)^{q}
\,dx\right)^{\frac{1}{q}}= C\left(\,\,\intl_{\Omega}
\|\nabla^kf\|^{q}\,dx\right)^{\frac{1}{q}}=
C\|f\|_{L^k_q(\Omega)}.
$$
Hence
$$
\|f^{\sharp}_{k,\Omega}\|_{L_q(\Omega)}\le
C(\|f\|_{L^k_q(\Omega)}+\|f\|_{L_q(\Omega)}) \le
C\|f\|_{\WKQO}.
$$
\par Theorem \reff{MAIN-EXT} is completely proved.\bx

\vspace*{10mm}
Department of Mathematics\\
Technion - Israel Institute of Technology\\
32000 Haifa\\
Israel\\\\
e-mail: pshv@tx.technion.ac.il\\\\

\begin{thebibliography}{ABCD}
\bibitem [AHHL]{AHHL} K. Astala, K. Hag, P. Hag, V.
    Lappalainen, Lipschitz Classes and the Hardy-Littlewood
    Property, Mh. Math 115 (1993), 267--279.
\bibitem [Br1]{Br1} Yu. A. Brudnyi, Spaces that are
    definable by means of local approximations, Trudy
    Moscov. Math. Obshch., 24 (1971) 69--132; English
    Transl.: Trans. Moscow Math Soc., 24 (1974) 73--139.
\bibitem [BrK]{BrK} Yu. A. Brudnyi, B.D. Kotljar, A certain
    problem of combinatorial geometry,  Sibirsk. Mat. Z. 11
    (1970) 1171–-1173; English transl. in Siberian Math. J.
    11 (1970), 870–-871.
\bibitem [BKos]{BKos} S. Buckley, P. Koskela, Criteria for
    imbeddings of Sobolev-Poincar\'{e} type, Internat.
    Math. Res. Notices 18 (1996), 881--902.
\bibitem [BSt2]{BSt2} S. Buckley, A. Stanoyevitch, Weak
    slice conditions and H\"older
    imbeddings, J. London Math. Soc. 66 (2001), 690--706.
\bibitem [BSt3]{BSt3} S. Buckley, A. Stanoyevitch, Weak
    slice conditions, product domains, and quasiconformal
    mappings, Rev. Math. Iberoam. 17 (2001), 1--37.
\bibitem [BSt]{BSt} S. Buckley, A. Stanoyevitch,
    Distinguishing properties of weak slice conditions,
    Conformal geometry and dynamics, An Electronic Journal
    of the AMS 7, (2003) 49-75.
\bibitem [BM]{BM} Yu. D. Burago, V. G. Maz'ya,  Certain
    Questions of Potential Theory and Function Theory for
    Regions with Irregular Boundaries. (Russian) Zap.
    Naucn. Sem. Leningrad. Otdel. Mat. Inst. Steklov.
    (LOMI) {\bf 3} 1967, 152 pp.; English trans.:
    Potential Theory and Function Theory for Irregular
    Regions. Seminars in Math., V. A. Steklov Math. Inst.,
    Leningrad, Vol. 3, Consultants Bureau, New York
    1969 vii+68 pp.
\bibitem [C2]{C2} A. P. Calder\'{o}n, Lebesgue spaces of
    differentiable functions and distributions, Proc. Symp.
    Pure Math., vol. IV, (1961) pp. 33--49.
\bibitem [C]{C} A. P. Calder\'{o}n, Estimates for singular
    integral operators in terms of maximal functions,
    Studia Math. {\bf 44,} (1972) 563--582 .
\bibitem [CF]{CF} R. Coifman and C. Fefferman, Weighted
    norm inequalities for maximal functions and singular
    integrals, Studia Math. 51 (1974) 241-–250.
\bibitem [GR]{GR} J. Garcia-Cuerva, J. L. Rubio de
    Francia, Weighted norm inequalities and related topics,
    North-Holland Mathematical Studies 116 (North-Holland,
    Amsterdam, 1985).
\bibitem [G]{G} F. W. Gehring, The $L^p$ integrability of
    the partial derivatives of a quasiconformal mapping,
    Acta Math. 130 (1973) 265–-277.
\bibitem [GM]{GM} F.W. Gehring, O. Martio, Lipschitz
    classes and quasiconformal mappings, Ann. Acad. Sci.
    Fenn. Ser. AI Math. 10 (1985), 203–-219.
\bibitem [GLV]{GLV} V.M. Gol'dshtein, T.G. Latfullin, S.K.
    Vodop'yanov, Criteria for extension of functions of the
    class $L^1_2$ from unbounded plain domains, Siber.
    Math. J.(Engl. transl.) 20, no. 2 (1979) 298-–301.
\bibitem [GV1]{GV1}  V. M. Gol'dstein , S. K.
    Vodop'janov,  Prolongement des fonctions de classe
    ${\cal L}\sp{1}\sb{p}$ et applications quasi
    conformes.,  C. R. Acad. Sci. Paris Ser. A-B 290
    (1980), no. 10, A453--A456.
\bibitem [GV2]{GV2}  V. M. Gol'dstein and S. K.
    Vodop'janov,  Prolongement des fonctions
    differentiables hors de domaines plans., C. R. Acad.
    Sci. Paris Ser. I Math. 293 (1981), no. 12, 581--584.
\bibitem [HKT]{HKT} P. Hajlasz, P. Koskela, H. Tuominen,
    Sobolev embeddings, extensions and measure density
    condition, J. Funct. Anal. 254 (2008) 1217-–1234.
\bibitem [HKT1]{HKT1} P. Hajlasz, P. Koskela, H. Tuominen,
    Measure density and extendability of Sobolev functions,
    Rev. Mat. Iberoamericana 24 (2008), no. 2, 645–-669.
\bibitem [Jn]{Jn} P. W. Jones, Quasiconformal mappings and
    extendability of functions in Sobolev spaces, Acta
    Math., 147 (1981), 71--78.
\bibitem [K]{K} P. Koskela, Extensions and Imbeddings,
    J. Funct. Anal. 159  (1998) 369--383.
\bibitem [L]{L}  V. Lappalainen, $Lip_h$-extension domains,
    Ann. Acad. Sci. Fenn. Ser. A I Math. Dissertations 56,
    (1985), 1-52.
\bibitem [Mac]{Mac} B. Muckenhoupt, Weighted norm
    inequalities
    for the Hardy–Littlewood maximal function, Trans.
    Amer. Math. Soc. 165 (1972) 207–-226.
\bibitem [M]{M}  V.G. Maz'ja,  Sobolev spaces,
     Springer-Verlag, Berlin, 1985, xix+486 pp.
\bibitem [MP]{MP} V. Maz'ya,  S. Poborchi, Differentiable
    Functions on Bad Domains, Word Scientific, River Edge,
    NJ, 1997.
\bibitem [Mel]{Mel} A. D. Melas, A sharp $L_p$ inequality
    for dyadic $A_1$ weights in $R^n$, Bull. London Math.
    Soc., 37 (2005) 919–-926.
\bibitem [S3]{S3} P. Shvartsman, On extension of Sobolev
    functions defined on regular subsets of metric measure
    spaces, J. Approx. Theory, 144 (2007), 139–-161.
\bibitem [S4]{S4} P. Shvartsman, Local approximations and
    intrinsic characterizations of spaces of smooth
    functions on regular subsets of $\RN$, Math. Nachr. 279
    (2006), no.11, 1212–-1241.
\bibitem [S7]{S7} P. Shvartsman, Sobolev $W^1_p$-spaces
    on closed subsets of $\RN$, Advances in Math.
    220 (2009) 1842–-1922.
\bibitem [St]{St} E. M. Stein, Singular integrals and
    differentiability properties of functions, Princeton
    Univ. Press, Princeton, New Jersey, 1970.
\bibitem [T2]{T2} H. Triebel, Theory of function spaces.
    II. Monographs in Mathematics, 84. Birkh\"{a}user
    Verlag, Basel, 1992.
\bibitem [W3]{W3} H. Whitney, Functions differentiable on
    the boundaries of regions, Ann. of Math. 35, no. 3,
    (1934), 482–-485.
\bibitem [Zob1]{Zob1} N. Zobin, Whitney's problem on
    extendability of functions and an intrinsic metric,
    {\it Advances in Math.} {\bf 133} (1998) 96--132.
\bibitem [Zob2]{Zob2} N. Zobin, Extension of smooth
    functions from finitely connected planar domains,
    J. Geom. Anal. 9, no. 3, (1999), 489--509.
\end{thebibliography}
\end{document}